\documentclass[twoside,11pt,a4paper,leqno]{article}
\usepackage{amssymb,amsmath,amscd,euscript,verbatim,array}
\usepackage[center,pagestyles]{titlesec}
\usepackage{geometry}
\usepackage{anysize}
\usepackage{fancyhdr}
\usepackage{indentfirst}
\usepackage{graphicx}
\usepackage{color}
\usepackage{ifpdf}
\usepackage{booktabs}
\usepackage{float}
\marginsize{3.5cm}{3.5cm}{2cm}{2cm}

\titleformat{\section}{\centering\normalsize}{\thesection.}{0.5em}{}
\titleformat{\subsection}{\normalsize\bfseries}{\thesubsection.}{0.5em}{}
\titleformat{\subsubsection}{\normalsize\bfseries}{\thesubsubsection.}{0.5em}{}
\titleformat{\paragraph}{\normalsize\bfseries}{\theparagraph.}{0.5em}{}
\titleformat{\subparagraph}{\normalsize\bfseries}{\thesubparagraph.}{0.5em}{}
 \titlespacing{\subparagraph}{-0.05em}{*1.5}{*1.1}
\newcommand{\N}{\mathbb{N}}
\newcommand{\Z}{\mathbb{Z}}
\newcommand{\R}{\mathbb{R}}
\newcommand{\Q}{\mathbb{Q}}

\newtheorem{Theorem}{Theorem}[section]
\newtheorem{Definition}[Theorem]{Definition}
\newtheorem{Lemma}[Theorem]{Lemma}
\newtheorem{Exercise}[Theorem]{Exercise}
\newtheorem{Proposition}[Theorem]{Proposition}
\newtheorem{Remark}[Theorem]{Remark}
\newtheorem{Corollary}[Theorem]{Corollary}

\newcommand{\eps}{\varepsilon}

\newcommand{\T}{\mathbb{T}}
\newcommand{\bthm}{\begin{Theorem}}
\newcommand{\ethm}{\end{Theorem}}
\newcommand{\bpr}{\begin{Proposition}}
\newcommand{\epr}{\end{Proposition}}
\newcommand{\blm}{\begin{Lemma}}
\newcommand{\elm}{\end{Lemma}}
\newcommand{\bex}{\begin{Exercise}}
\newcommand{\eex}{\end{Exercise}}
\newcommand{\be}{\begin{equation}}
\newcommand{\ee}{\end{equation}}
\newcommand{\beal}{\begin{aligned}}
\newcommand{\enal}{\end{aligned}}
\newcommand{\brm}{\begin{Remark}}
\newcommand{\erm}{\end{Remark}}
\newcounter{item}[section]

\newcommand{\Proof}{\textbf{Proof}\hspace{0.3cm}}
\newcommand{\End}{\ensuremath{\hfill{\Box}}\\}
\renewcommand{\title}[1]{\begin{center}\textbf{\large #1}\end{center}}
\renewcommand{\author}[1]{\begin{center}\small #1\end{center}}
\renewcommand{\date}[1]{\begin{center}#1\end{center}}

\makeatletter \@addtoreset{equation}{subsection}
\makeatother
\begin{document}
\vspace{10pt}
\title{CONVERSE KAM THEORY REVISITED}
\vspace{6pt}
\author{\sc Lin Wang}
\vspace{10pt} \thispagestyle{plain}
{\begingroup\makeatletter
\let\@makefnmark\relax  \footnotetext{{\it Date}: \today}
\makeatother\endgroup}
\begin{quote}
\small {\sc Abstract.} For an integrable Hamiltonian with $d\ (d\geq
2)$ degrees of freedom, we show the conditions on perturbations, for
which invariant tori can be destructed.
\end{quote}

\tableofcontents
\section{\sc Introduction}
By the Kolmogorov, Arnold and Moser (KAM) theory, we know that under
certain non-degeneracy, most (full Lebesgue measure) invariant tori
of an integrable Hamiltonian system are persisted under small
perturbations. As the sizes of the perturbations increase, those
persisted invariant tori are destructed progressively. The problems
of determining the critical boundary between the persistence and
destruction of invariant tori motivates so called converse KAM
theory.

Roughly speaking, the converse KAM theory consists of two parts. The
first part, with more physical flavor, is concerned about
destruction of invariant tori under the perturbations with positive
lower bound, which is started by the analysis of invariant circles
of the generalized standard maps (see e.g. \cite{Chi},\cite{H1} and
\cite{Mm2}). Besides, numerical results go further than theoretical
ones. More sharp boundaries are obtained by numerical methods (see
e.g. \cite{Gr} and \cite{MP}). The second part is concerned about
destruction of invariant tori under arbitrarily small perturbations
in certain topology and it seems that numerical method is not as
efficient as the first part. In this paper, we are devoted to
develop the second part.

In 1962, Moser proved that the invariant circles with Diophantine
rotation numbers of an integrable twist map is persisted under
arbitrarily small perturbations in the $C^{333}$ topology
(\cite{Mo1}). By the efforts of Moser, R\"{u}ssman, Herman and
P\"{o}schel (\cite{H1,H33},\cite{Mo2,Mo3},\cite{P} and
\cite{R1,R2}), for Hamiltonian systems with $d$-degrees of freedom,
it is obtained that certain invariant tori are persisted under
arbitrarily small perturbations in the $C^{2d+\delta}$ topology,
where $\delta$ is a small positive constant. Especially, Herman
proved in \cite{H33} that for twist maps on annulus, certain
invariant circles can be persisted under arbitrarily small
perturbations in the $C^3$ topology.

In contrast with the results on persistence of invariant tori, for
exact area-preserving twist maps on annulus, it is proved by Herman
in \cite{H1} that invariant circles with given rotation numbers can
be destructed by $C^{3-\delta}$ arbitrarily small $C^\infty$
perturbations.  For certain rotation numbers, it is obtained by
Mather (resp. Forni) in \cite{M4} (resp. \cite{Fo}) that the
invariant circles with those rotation numbers can be destroyed by
small perturbations in finer topology respectively. More precisely,
Mather considers Liouville rotation numbers and the topology of the
perturbation induced by $C^\infty$ metric. Forni is concerned about
more special rotation numbers which can be approximated  by rational
ones exponentially and the topology of the perturbation induced by
the supremum norm of $C^\omega$ (real-analytic) function.  Bessi
extended the result to the systems with multi-degrees of freedom. He
found that the invariant Lagrangian torus with certain rotation
vector can be destructed by an arbitrarily small $C^\omega$
perturbation for certain positive definite systems with
multi-degrees of freedom in \cite{B2}, where Lagrangian torus is a
natural analogy to invariant circle in multi-degrees of freedom (see
Definition \ref{dd1} below).

On the other hand, it is also proved by Herman in \cite{H3} that all
Lagrangian tori of an integrable symplectic twist map with $d\geq 1$
degrees of freedom can be destructed by $C^{d+2-\delta}$ arbitrarily
small $C^\infty$ perturbations of the generating function.
Equivalently (\cite{Go,Mo4}), it shows that all Lagrangian tori of
an integrable Hamiltonian system with $d\geq 2$ degrees of freedom
can be destructed by $C^\infty$ perturbations which are arbitrarily
small in the $C^{d+1-\delta}$ topology.  Roughly speaking, there is
a balance among the arithmetic property of the rotation vector, the
regularity of the perturbation and its topology.

Comparing the results on both sides, it is natural to ask the
following questions: \begin{itemize}
\item for every given rotation vector $\omega$,
if the Lagrangian torus with $\omega$ can be destructed by an
arbitrarily small $C^\infty$ perturbation in the $C^r$ topology,
then what is the maximum of $r$?
\item for every given rotation vector $\omega$,
if the Lagrangian torus with $\omega$ can be destructed by an
arbitrarily small analytic perturbation in the $C^r$ topology, then
what is the maximum of $r$?
\item  if all of Lagrangian tori can
be destructed by an arbitrarily small real-analytic perturbation in
the $C^r$ topology, then what is the maximum of $r$?
\end{itemize}

Based on \cite{CW} and \cite{W1,W2,W3,W4}, we have the following
theorems.
\begin{Theorem}\label{T1} Given an integrable Hamiltonian $H_0$ with $d\ (d\geq 2)$
 degrees of freedom and a
rotation vector $\omega$, there exists a sequence of $C^\infty$
Hamiltonians $\{H_n\}_{n\in \N}$ such that $H_n\rightarrow H_0$ in
$C^{2d-\delta}$ topology and
 the Hamiltonian flow generated by $H_n$ does not admit the
Lagrangian torus with the rotation vector $\omega$.
\end{Theorem}
\begin{Theorem}\label{T2} Given an integrable Hamiltonian $H_0$ with $d\ (d\geq 2)$
 degrees of freedom and a
rotation vector $\omega$, there exists a sequence of $C^\omega$
Hamiltonians $\{H_n\}_{n\in \N}$ such that $H_n\rightarrow H_0$ in
$C^{d+1-\delta}$ topology and
 the Hamiltonian flow generated by $H_n$ does not admit the
Lagrangian torus with the rotation vector $\omega$.
\end{Theorem}
\begin{Theorem}\label{T3}
 For an integrable Hamiltonian $H_0$ with $d\ (d\geq 2)$ degrees
of freedom, all Lagrangian tori can be destructed by analytic
perturbations which are arbitrarily small in the $C^{d-\delta}$
topology.
\end{Theorem}
For positive definite Hamiltonian systems with $d\geq 2$ degrees of
freedom, together with Herman's result in \cite{H3}, we have the
following Table 1.
\begin{table}[!h]
\centering
\begin{tabular}{c|c|c}
  $r$& $C^\infty$ & $C^\omega$ \\
   \hline
  Single Destruction & $2d-\delta$ & $d+1-\delta$ \\
   \hline
  Total Destruction & $d+1-\delta$ & $d-\delta$ \\
  \hline
\end{tabular}
\caption{Values of $r$ for destruction of Lagrangian torus (tori) in
the $C^r$ topology}
\end{table}

Unfortunately, except for the destruction of a Lagrangian torus by
the $C^\infty$ perturbations,
 we still don't know
whether the other results are optimal. Some further developments of
KAM theory are needed to verify the optimality.

This paper is outlined as follows. In Section 2, we consider the
destruction of a Lagrangian torus with given rotation vector
(Theorem \ref{T1} and Theorem \ref{T2}). Based on the difference of
topology of phase spaces between $d=2$ and $d\geq 3$, we divide the
arguments into two cases. In terms of the correspondence between
exact area-preserving twist maps and Hamiltonian systems with 2
degrees of freedom, the problem on destruction of the Lagrangian
torus for the Hamiltonian system is transformed into the one on
destruction of the invariant circle for the twist map. Using
variational method developed by Mather, a new proof of Herman's
result (\cite{H1}) is provided. Moreover, from Jackson's
approximation, the result on $C^\omega$ perturbation is obtained.
For the case with $d\geq 3$ degrees of freedom, the minimality of
the orbits on Lagrangian torus plays a crucial role in destruction
of the torus. Combining with Melnikov method, the destruction of
Lagrangian torus under $C^\omega$ perturbations is achieved. In
Section 3, we are concerned about the destruction of all Lagrangian
tori (Theorem \ref{T3}). Similar to Section 2, according to the
correspondence between the Hamiltonian system and the symplectic
twist map, we focus on destruction of Lagrangian tori for symplectic
twist map. From a criterion of total destruction of Lagrangian tori
found by Herman and an approximation lemma, total destruction of
Lagrangian tori under $C^\omega$ perturbations is completed.
\section{\sc Destruction of a Lagrangian torus}

\subsection{Case with $2$ degrees of freedom}
Based on the correspondence between exact area-preserving twist maps
and Hamiltonian systems with 2 degrees of freedom, it is sufficient
to consider the destruction of invariant circle for exact
area-preserving twist map.
 \subsubsection{ Preliminaries of exact
area-preserving twist map}

\paragraph{Minimal configuration} Let $f$: $\T\times\R\rightarrow
\T\times\R$ ($\T=\R/\Z$) be an exact area-preserving monotone twist
map and $h$: $\R^2\rightarrow\R^2$ be a generating function for the
lift $F$ of $f$ to $\R^2$, namely $F$ is generated by the following
equations
\begin{equation*}
\begin{cases}
y=-\partial_1 h(x,x'),\\
y'=\partial_2 h(x,x'),
\end{cases}
\end{equation*}
where $F(x,y)=(x',y')$. The lift $F$ gives rise to a dynamical
system whose orbits are given by the images of points of $\R^2$
under the successive iterates of $F$. The orbit of the point
$(x_0,y_0)$ is the bi-infinite sequence
\[\{...,(x_{-k},y_{-k}),...,(x_{-1},y_{-1}),(x_0,y_0),(x_1,y_1),...,(x_k,y_k),...\},\]
where $(x_k,y_k)=F(x_{k-1},y_{k-1})$. The sequence
\[(...,x_{-k},...,x_{-1},x_0,x_1,...,x_k,...)\] denoted by $(x_i)_{i\in\Z}$ is called
stationary configuration if it stratifies the identity
\[\partial_1 h(x_i,x_{i+1})+\partial_2 h(x_{i-1},x_i)=0,\ \text{for\ every\ }i\in\Z.\]
Given a sequence of points $(z_i,...,z_j)$, we can associate its
action
\[h(z_i,...,z_j)=\sum_{i\leq s<j}h(z_s,z_{s+1}).\] A configuration $(x_i)_{i\in\Z}$
is called minimal if for any $i<j\in \Z$, the segment of
$(x_i,...,x_j)$ minimizes $h(z_i,...,z_j)$ among all segments
$(z_i,...,z_j)$ of the configuration  satisfying $z_i=x_i$ and
$z_j=x_j$. It is easy to see that every minimal configuration is a
stationary configuration. By \cite{B}, minimal configurations
satisfy a group of remarkable properties as follows:
\begin{itemize}
\item Two distinct minimal configurations cross at most once, which
is so called Aubry's crossing lemma.
\item For every minimal configuration $\bold{x}=(x_i)_{i\in\Z}$, the limit
\[\rho(\bold{x})=\lim_{n\rightarrow\infty}\frac{x_{i+n}-x_i}{n}\]
exists and doesn't depend on $i\in\Z$. $\rho(\bold{x})$ is called
the rotation number of $\bold{x}$.
\item For every $\omega\in \R$, there exists a minimal configuration
with rotation number $\omega$. Following the notations of [B], the
set of all minimal configurations with rotation number $\omega$ is
denoted by $M_\omega^h$, which can be endowed with the topology
induced from the product topology on $\R^\Z$. If
$\bold{x}=(x_i)_{i\in\Z}$ is a minimal configuration, considering
the projection $pr:\ M_\omega^h\rightarrow\R$ defined by
$pr(\bold{x})=x_0$, we set $\mathcal {A}_\omega^h=pr(M_\omega^h)$.
\item If $\omega\in\Q$, say $\omega=p/q$ (in lowest terms), then it is convenient to define the rotation symbol to detect the structure of
$M_{p/q}^h$. If $\bold{x}$ is a minimal configuration with rotation
number $p/q $, then the rotation symbol $\sigma(\bold{x})$ of
$\bold{x}$ is defined as follows
\begin{equation*}
\sigma(\bold{x})=\left\{\begin{array}{ll}
\hspace{-0.4em}p/q+,&\text{if}\ x_{i+q}>x_i+p\ \text{for\ all\ }i,\\
\hspace{-0.4em}p/q,&\text{if}\ x_{i+q}=x_i+p\ \text{for\ all\ }i,\\
\hspace{-0.4em}p/q-,&\text{if}\ x_{i+q}<x_i+p\ \text{for\ all\ }i.\\
\end{array}\right.
\end{equation*}
 Moreover, we set
\begin{align*}
&M_{{p/q}^+}^h=\{\bold{x} \text{\ a is minimal configuration with
rotation symbol}\  p/q \text{\ or\ } p/q+\},\\
&M_{{p/q}^-}^h=\{\bold{x} \text{\ a is minimal configuration with
rotation symbol}\  p/q \text{\ or\ } p/q-\},
\end{align*}
then both $M_{{p/q}^+}^h$ and $M_{{p/q}^+}^h$ are totally ordered.
Namely, every two configurations in each of them do not cross. We
denote $pr(M_{{p/q}^+}^h)$ and $pr(M_{{p/q}^-}^h)$ by $\mathcal
{A}_{{p/q}^+}^h$ and $\mathcal {A}_{{p/q}^-}^h$ respectively.
\item If $\omega\in\R\backslash\Q$ and $\bold{x}$ is a minimal
configuration with rotation number $\omega$, then
$\sigma(\bold{x})=\omega$ and $M_\omega^h$ is totally ordered.
\item $\mathcal {A}_\omega^h$ is a closed subset of $\R$ for every rotation symbol
$\omega$.
\end{itemize}
\paragraph{Peierls's barrier} In \cite{M3}, Mather introduced the notion of
Peierls's barrier and gave a criterion of existence of invariant
circle. Namely, the exact area-preserving monotone twist map
generated by $h$ admits an invariant circle with rotation number
$\omega$ if and only if the Peierls's barrier $P_\omega^h(\xi)$
vanishes identically for all $\xi\in\R$. The Peierls's barrier is
defined as follows:
\begin{itemize}
\item If $\xi\in \mathcal {A}_\omega^h$, we set $P_\omega^h(\xi)$=0.
\item If $\xi \not\in \mathcal {A}_\omega^h$, since $\mathcal {A}_\omega^h$ is a closed set in $\R$, then $\xi$ is contained in some
complementary interval $(\xi^-,\xi^+)$ of $\mathcal {A}_\omega^h$ in
$\R$. By the definition of $\mathcal {A}_\omega^h$, there exist
minimal configurations with rotation symbol $\omega$,
$\bold{x^-}=(x_i^-)_{i\in\Z}$ and $\bold{x^+}=(x_i^+)_{i\in\Z}$
satisfying $x_0^-=\xi^-$ and $x_0^+=\xi^+$. For every configuration
$\bold{x}=(x_i)_{i\in\Z}$ satisfying $x_i^-\leq x_i\leq x_i^+$, we
set
\[G_\omega(\bold{x})=\sum_I(h(x_i,x_{i+1})-h(x_i^-,x_{i+1}^-)),\]
where $I=\Z$, if $\omega$ is not a rational number, and $I=\{0,...,
q-1 \}$, if $\omega=p/q$. $P_\omega^h(\xi)$ is defined as the
minimum of $G_\omega(\bold{x})$ over the configurations $\bold{x}\in
\Pi=\prod_{i\in I}[x_i^-,x_i^+]$ satisfying $x_0=\xi$. Namely
\[P_\omega^h(\xi)=\min_{\bold{x}}\{G_\omega(\bold{x})|\bold{x}\in \Pi\ \text{and}\ \ x_0=\xi\}.\]
\end{itemize}
By \cite{M3}, $P_\omega^h(\xi)$ is a non-negative periodic function
of the variable $\xi\in\R$ with the modulus of continuity with
respect to $\omega$.

\subsubsection{$C^\infty$ case}
\paragraph{Construction of the generating functions}

Consider a completely integrable system with the generating function
\[h_0(x,x')=\frac{1}{2}(x-x')^2 \quad x,x'\in \R.\] We construct the perturbation consisting of two
parts. The first one is
\begin{equation}\label{31}
u_n(x)=\frac{1}{n^a}(1-\cos(2\pi x) )\quad x\in \R,\end{equation}
where $n\in \N$ and $a$ is a positive constant independent of $n$.
The second one is a non negative function $v_n(x)$ satisfying
\begin{equation}\label{32}
\begin{cases}
 v_n(x+1)=v_n(x),\\
\text{supp}\,v_n\cap [0,1]\subset
[\frac{1}{2}-\frac{1}{n^a},\frac{1}{2}+\frac{1}{n^a}], \\
\max v_n=n^{-s}, \\
{||v_n||}_{C^k}\sim n^{-s'},
\end{cases}
\end{equation}
where  $f\sim g$ means that $\frac{1}{C}g<f<Cg$ holds for a constant
$C>1$. For further deduction, we need $s'>a$. It is enough to take
$s=(k+2)a$ for achieving that. The generating function of the nearly
integrable system is constructed as follow:
\begin{equation}\label{h}
h_n(x,x')=h_0(x,x')+u_n(x')+v_n(x'),
\end{equation}
where $n\in\N$. Moreover, we have the following theorem.

\bthm\label{MR} For $\omega\in\R\backslash\Q$ and $n$ large enough,
the exact area-preserving monotone twist map generated by $h_n$ does
not admit any invariant circles with the rotation number satisfying
\[|\omega|<n^{-\frac{a}{2}-\delta}, \] where $\delta$ is a small positive constant independent of $n$. \ethm
We will prove Theorem \ref{MR} in the following sections. First of
all, based on the theorem, we verify that our example has the
property aforementioned in Section 1.

If $\omega\in \Q$, then the invariant circles with rotation number
$\omega$ could be easily destructed by an analytic perturbation
arbitrarily close to $0$. Therefore it suffices to consider the
irrational $\omega$. The case with a given irrational rotation
number can be easily reduced to the one with a small enough rotation
number. More precisely,

\blm\label{Herm} Let $h_P$ be a generating function as follow
\[h_P(x,x')=h_0(x,x')+P(x'),\] where $P$ is a periodic
function of periodic $1$. Let $Q(x)=q^{-2}P(qx),q\in \N$, then the
exact area-preserving monotone twist map generated by
$h_Q(x,x')=h_0(x,x')+Q(x')$ admits an invariant circle with rotation
number $\omega \in \R\backslash \Q$ if and only if the exact
area-preserving monotone twist map generated by $h_P$ admits an
invariant circle with rotation number $q\omega-p, p\in \Z$. \elm

We omit the proof and for more details, see [H2]. For the sake of
simplicity of notations, we denote $Q_{q_n}$ by $Q_n$ and the same
to $u_{q_n}, v_{q_n}$ and $h_{q_n}$. Let
\[Q_n(x)={q_n}^{-2}(u_n(q_nx)+v_n(q_nx)),\] where $(q_n)_{n\in \N}$
is a sequence satisfying Dirichlet approximation
\begin{equation}\label{diri}
|q_n\omega-p_n|<\frac{1}{q_n}, \end{equation} where $p_n\in \Z$ and
$q_n\in \N$. Since $\omega\in\R\backslash\Q$, we say
$q_n\rightarrow\infty$ as $n\rightarrow\infty$. Let
$\tilde{h}_n(x,x')=h_0(x,x')+Q_n(x')$, we have

\begin{Corollary}\label{Mcor}
For a given rotation number $\omega\in \R\backslash\Q$ and every
$\eps$, there exists $N$ such that for $n>N$, the exact
area-preserving monotone map generated by $\tilde{h}_n$ admits no
invariant circle with rotation number $\omega$ and
\[||\tilde{h}_n-h_0||_{C^{4-\delta'}}<\eps,\]where $\delta'$ is a small positive constant independent of $n$.
\end{Corollary}

\Proof Based on Theorem \ref{MR} and Dirichlet approximation
(\ref{diri}), it suffices to take
\[\frac{1}{q_n}\leq\frac{1}{{q_n}^{\frac{a}{2}+\delta}},\]
which implies \begin{equation}\label{33}
a\leq2-2\delta.\end{equation} From (\ref{31}), (\ref{32}) and
(\ref{h}), it follows that
\begin{align*}
 ||\tilde{h}_n&(x,x')-h_0(x,x')||_{C^r}\\
 &=||Q_n(x')||_{C^r},\\
&\leq{q_n}^{-2}(||u_n(q_nx')||_{C^r}+||v_n(q_nx')||_{C^r}),\\
&\leq{q_n}^{-2}({q_n}^{-a}(2\pi)^r{q_n}^r+C_1{q_n}^{-s'}{q_n}^r),\\
&\leq C_2{q_n}^{r-a-2},
\end{align*}
where $C_1, C_2$ are positive constants only depending on $r$.

To complete the proof, it is enough to make $r-a-2<0$, which
together with (\ref{33}) implies
\[r<a+2\leq 4-2\delta.\] We set ${\delta}'=2\delta$, then the
proof of Corollary \ref{Mcor} is completed.\End

The following sections are devoted to prove Theorem \ref{MR}. For
simplicity, we don't distinguish the constant $C$ in following
different estimate formulas.

\paragraph{Estimate of lower bound of $P_{0^+}^{h_n}$ }

In this section, we will estimate the lower bound of $P_{0^+}^{h_n}$
at a given point. To achieve that, we need to estimate the distances
of pairwise adjacent elements of the minimal configuration.
\subparagraph{A spacing lemma}

\blm\label{lowstep} Let $(x_i)_{i\in \Z}$ be a minimal configuration
of $\bar{h}_n$ with rotation symbol $0^+$, then
\[x_{i+1}-x_i\geq C(n^{-\frac{a}{2}}),\quad \text{for}\quad x_i\in \left[\frac{1}{4},\frac{3}{4}\right],\] where
$\bar{h}_n(x_i,x_{i+1})=h_0(x_i,x_{i+1})+u_n(x_{i+1})$. \elm

\Proof Without loss of generality, we assume $x_i\in [0,1]$ for all
$i\in\Z$. By Aubry's crossing lemma, we have
\[0<...<x_{i-1}<x_i<x_{i+1}<...<1.\]We consider the configuration
$(\xi_i)_{i\in \Z}$ defined by
\begin{equation*}
\xi_j= \left\{\begin{array}{ll}\hspace{-0.4em}x_j,& j<i,\\
\hspace{-0.4em}x_{j+1},& j\geq i.\\
\end{array}\right.
\end{equation*}
Since $(x_i)_{i\in \Z}$ is minimal, we have
\[\sum_{i\in \Z}\bar{h}_n(\xi_i,\xi_{i+1})-\sum_{i\in \Z}\bar{h}_n(x_i,x_{i+1})\geq 0.\]
By the definitions of $\bar{h}_n$ and $(\xi_i)_{i\in\Z}$, we have
\begin{align*}
0&\leq\sum_{i\in \Z}\bar{h}_n(\xi_i,\xi_{i+1})-\sum_{i\in
\Z}\bar{h}_n(x_i,x_{i+1})\\
&=\bar{h}_n(x_{i-1},x_{i+1})-\bar{h}_n(x_{i-1},x_{i})-\bar{h}_n(x_{i},x_{i+1})\\
&=(x_{i+1}-x_i)(x_i-x_{i-1})-u_n(x_i).
\end{align*}
Moreover,\[u_n(x_i)\leq(x_{i+1}-x_i)(x_i-x_{i-1})\leq\frac{1}{4}(x_{i+1}-x_{i-1})^2.\]
Therefore, \[x_{i+1}-x_{i-1}\geq 2\sqrt{u_n(x_i)}.\] For $x_i\in
[\frac{1}{4},\frac{3}{4}]$, $u_n(x_i)\geq n^{-a}$, hence,
\begin{equation}\label{ls} x_{i+1}-x_{i-1}\geq 2n^{-\frac{a}{2}}.
\end{equation}

Since $(x_i)_{i\in\Z}$ is a stationary configuration, we have
\begin{align*}
x_{i+1}-x_i&=-\partial_1\bar{h}_n(x_i,x_{i+1}),\\
&=\partial_2\bar{h}_n(x_{i-1},x_i),\\
&=x_i-x_{i-1}+u_n'(x_i).
\end{align*}
Since $u_n'(x)=\frac{2\pi}{n^a}\sin(2\pi x)$, it follows from
$(\ref{ls})$ that
\[x_{i+1}-x_i\geq C(n^{-\frac{a}{2}}),\quad x_i\in \left[\frac{1}{4},\frac{3}{4}\right]
.\] The proof of Lemma \ref{lowstep} is completed.\End
\subparagraph{The lower bound of $P_{0^+}^{h_n}$ } By the definition
of $v_n$,
$\text{supp}\,v_n\cap[0,1]\subset[\frac{1}{2}-\frac{1}{n^a},\frac{1}{2}+\frac{1}{n^a}]$
and $v_n(x+1)=v_n(x)$. Let $(x_i)_{i\in \Z}$ be the minimal
configuration of
$\bar{h}_n(x_i,x_{i+1})=h_0(x_i,x_{i+1})+u_n(x_{i+1})$ with rotation
symbol $0^+$ satisfying $x_0=\frac{1}{2}-\frac{1}{n^a}$, then
\[(x_i)_{i\in \Z}\cap\text{supp}v_n=\emptyset.\]
Moreover, for all $i\in\Z$, \[v_n(x_i)=0.\]

 Let $(\xi_i)_{i\in \Z}$
be a minimal configuration of $h_n$ defined by (\ref{h}) with
rotation symbol $0^+$ satisfying $\xi_0=\eta$, where $\eta$
satisfies $v_n(\eta)=\max v_n(x)=n^{-s}$, then
\begin{align*}
\sum_{i\in
\Z}(h_n(&\xi_i,\xi_{i+1})-h_n(\xi_i^-,\xi_{i+1}^-))\\
&\geq v_n(\eta)+\sum_{i\in \Z}\bar{h}_n(\xi_i,\xi_{i+1})-\sum_{i\in
\Z}h_n(\xi_i^-,\xi_{i+1}^-),\\
&\geq v_n(\eta)+\sum_{i\in \Z}\bar{h}_n(x_i,x_{i+1})-\sum_{i\in
\Z}h_n(x_i,x_{i+1}),\\
&=v_n(\eta)-\sum_{i\in \Z}v_n(x_{i+1}),\\
&=v_n(\eta).
\end{align*}
Therefore,
\[P_{0^+}^{h_n}(\eta)=\min_{x_0=\eta}\sum_{i\in
\Z}(h_n(x_i,x_{i+1})-h_n(x_i^-,x_{i+1}^-))\geq v_n(\eta)=n^{-s}.\]We
conclude that there exists a point $\xi\in
[\frac{1}{2}-\frac{1}{n^a},\frac{1}{2}+\frac{1}{n^a}]$ such that
\begin{equation}\label{lowbound}
P_{0^+}^{h_n}(\xi)\geq n^{-s}.
\end{equation}

\paragraph{The approximation from $P_{0^+}^{h_n}$ to
$P_{\omega}^{h_n}$ } In this section, we will prove the improvement
of modulus of continuity of Peierls's barrier based on the
hyperbolicity of $h_n$. Namely

\blm\label{apprx} For every irrational rotation symbol $\omega$
satisfying $0<\omega<n^{-\frac{a}{2}-\delta}$, we have
\[|P_{\omega}^{h_n}(\xi)-P_{0^+}^{h_n}(\xi)|\leq C\exp\left(-2n^{\frac{\delta}{2}}\right).\]
where $\xi\in
\left[\frac{1}{2}-\frac{1}{n^a},\frac{1}{2}+\frac{1}{n^a}\right]$
and $\delta$ is a small positive constant independent of  $n$.
 \elm

\subparagraph{Some counting lemmas} To prove the lemma, we need to
do some preliminary work. First of all, we count the number of the
elements of a minimal configuration $(x_i)_{i\in\Z}$ with arbitrary
rotation symbol $\omega$ in a given interval. With the method of
[F], we can conclude the following lemma.

\blm\label{count 0} Let $(x_i)_{i\in\Z}$ be a minimal configuration
of $h_n$ with rotation symbol $\omega>0$,
$J_n=\left[\exp\left(-n^{\frac{\delta}{2}}\right),\frac{1}{2}\right]$
and $\Lambda_n=\{i\in \Z|\,x_i\in J_n\}$, then \[\sharp\Lambda_n\leq
Cn^{\frac{a}{2}+\frac{\delta}{2}},\] where $\sharp\Lambda_n$ denotes
the number of elements in $\Lambda_n$ and  $\delta$ is a small
positive constant independent of $n$. \elm

\Proof Let $x^-=\exp\left(-n^{\frac{\delta}{2}}\right),
x^+=\frac{1}{2}$ and
$\sigma=\left(\frac{x^+}{x^-}\right)^{\frac{1}{N}}$, hence,
\[\ln \sigma=\frac{\ln(x^+)-\ln(x^-)}{N}.\] We choose $N\in \N$
such that $1\leq\ln\sigma\leq 2$, then
$N=\Omega\left(n^{\frac{\delta}{2}}\right)$.

We consider the partition of the interval $J_n=[x^-,x^+]$ into the
subintervals $J_n^k=[\sigma^k x^-,\sigma^{k+1}x^-]$ where $0\leq
k<N$. Hence, $J_n=\cup_{k=0}^{N-1}J_n^k$. We set
$S_k=\{i\in\Lambda_n|(x_{i-1},x_{i+1})\subset J_n^k\}$ and
$m_k=\sharp S_k$.

By the similar deduction as the one in Lemma \ref{lowstep}, we have
\[x_{i+1}-x_{i-1}\geq 2\sqrt{u_n(x_i)+v_n(x_i)}\geq Cn^{-\frac{a}{2}}x_i,\quad\text{for}\quad x_i\in \left[0,\frac{1}{2}\right].\]
For simplicity of notation, we denote $Cn^{-\frac{a}{2}}$ by
$\alpha_n$.

If  there exists $k$ such that $i\in S_k$ for $(x_i)_{i\in \Z}$,
then $x_{i+1}-x_{i-1}\geq\alpha_n\sigma^kx^-$,
moreover,\[m_k\alpha_n\sigma^kx^-\leq 2\mathcal {L}(J_n^k)=
2(\sigma-1)\sigma^kx^-,\]where $\mathcal {L}(J_n^k)$ denotes the
length of the interval of $J_n^k$. Hence $m_k\leq 2(\sigma
-1)\alpha_n^{-1}$.

On the other hand, if $i\in \Lambda_n\backslash \cup_{k=0}^{N-1}S_k$
, then there exists $l$ satisfying $0\leq l<N$ such that
\[x_{i-1}<\sigma^lx^-<x_{i+1}.\]Hence,
\[\sharp\{i\in\Lambda_n|i\not\in S_k \ \text{for\ any}\ k\}\leq 2N.
\] Therefore,
\begin{align*}
\sharp(\Lambda_n)\leq 2N(\sigma-1)\alpha_n^{-1}+2N.
\end{align*}Since $1\leq\ln\sigma\leq 2$ and $N=\Omega\left(n^{\frac{\delta}{2}}\right)$, then we have
\[\sharp\Lambda_n\leq
Cn^{\frac{a}{2}+\frac{\delta}{2}}.\] The proof of Lemma \ref{count
0} is completed.\End
\begin{Remark}\label{rr1}
Let $(x_i)_{i\in\Z}$ be a minimal configuration of $h_n$ defined by
$(\ref{h})$ with rotation symbol $\omega>0$, An argument as similar
as the one in Lemma \ref{count 0} implies that
\[\sharp\left\{i\in\Z\bigg|\,x_i\in\left[\exp\left(-n^{\frac{\delta}{2}}\right),1-\exp\left(-n^{\frac{\delta}{2}}\right)\right]\right\}\leq Cn^{\frac{a}{2}+\frac{\delta}{2}}.\]
\end{Remark}

It is easy to count the number of the elements of a minimal
configuration with irrational rotation symbol. More precisely, we
have the following lemma.
 \blm\label{count w}Let
$(x_i)_{i\in \Z}$ be a minimal configuration with rotation number
$\omega\in \R\backslash\Q$. Then for every interval $I_k$ of length
$k$, $k\in \N$,
\[\frac{k}{\omega}-1\leq\sharp\{i\in \Z|x_i\in
I_k\}\leq\frac{k}{\omega}+1.\] \elm

\Proof For every minimal configuration $(x_i)_{i\in\Z}$ with
rotation number $\omega$, there exists an orientation-preserving
circle homeomorphism $\phi$ such that $\rho(\Phi)=\omega$, where
$\Phi:\R\rightarrow\R$ denotes a lift of $\phi$. Since $\omega\in
\R\backslash\Q$, thanks to [H1], $\phi$ has a unique invariant
probability measure $\bar{\mu}$ on $\T$ such that $\int_x^{\Phi(x)}
d\bar{\mu}=\omega$ for every $x\in\R$. We denote $\int_x^{\Phi(x)}
d\bar{\mu}$ by $\mu(x,\Phi(x))$. In particular,
\[\mu(x_i,x_{i+1})=\omega,\quad \text{for\ every\ }i\in\Z.\]
From $\mu(I_k)=k$, it follow that
\begin{align*}
&\omega(\sharp\{i\in \Z|x_i\in I_k\}-1)\leq k,\\
&\omega(\sharp\{i\in \Z|x_i\in I_k\}+1)\geq k,
 \end{align*}
which completes the proof of Lemma \ref{count w}.\End

Based on Lemma \ref{count 0} and Lemma \ref{count w}, if
$0<\omega<n^{-\frac{a}{2}-\delta}$ and $\omega$ is irrational, then
\begin{equation}\label{pnew}
\sharp\{i\in \Z|x_i\in I_1\}\geq \frac{1}{\omega}-1\geq
C_1n^{\frac{a}{2}+\delta}>C_2n^{\frac{a}{2}+\frac{\delta}{2}},
\end{equation}
where $I_1$ denotes the closed interval of length $1$. Moreover, we
have the following conclusion.

\blm\label{lbspace} Let $(x_i)_{i\in \Z}$ be a minimal configuration
of $h_n$ defined by $(\ref{h})$ with rotation symbol
$0<\omega<n^{-\frac{a}{2}-\delta}$, then there exist $j^-,j^+\in\Z$
such that \begin{align*}&0<x_{j^--1}<x_{j^-}<x_{j^-+1}\leq
\exp(-n^{\frac{\delta}{2}}),\\
& 1-\exp(-n^{\frac{\delta}{2}})\leq
x_{j^+-1}<x_{j^+}<x_{j^++1}<1.\end{align*}\elm

\Proof By contradiction, we assume that there exist at most two
points of $(x_i)_{i\in\Z}$ in $[0, \exp(-n^{\frac{\delta}{2}})]$,
say $x_m$ and $x_{m+1}$. It follows that $x_{m-1}<0$ and
$x_{m+2}>\exp(-n^{\frac{\delta}{2}})$. Hence, among the intervals
$[x_{m-1}, x_m]$, $[x_m, x_{m+1}]$ and $[x_{m+1}, x_{m+2}]$, there
exists at least one such that its length is not less than
$\frac{1}{3}\exp(-n^{\frac{\delta}{2}})$. Without loss of
generality, say $[x_{m+1}, x_{m+2}]$.

Since $(x_i)_{i\in\Z}$ is a stationary configuration, we have
\[x_{m+2}-x_{m+1}=x_{m+1}-x_m+u'_n(x_{m+1}),\] where $u'_n(x_{m+1})=\frac{2\pi}{n^a}\sin(2\pi
x_{m+1})$. From $x_{m+1}\in [-\exp(-n^{\frac{\delta}{2}}),
\exp(-n^{\frac{\delta}{2}})]$, it follows that
\[|u'_n(x_{m+1})|\leq
Cn^{-a}\exp(-n^{\frac{\delta}{2}}),\]which implies there exists $N$
independent of $n$ such that $[-\exp(-n^{\frac{\delta}{2}}),
\exp(-n^{\frac{\delta}{2}})]$ contains at most $N$ points of
$(x_i)_{i\in \Z}$.

On the other hand, by (\ref{pnew}), we have that for $n$ large
enough, the number of points of $(x_i)_{i\in\Z}$ in
$[-\exp(-n^{\frac{\delta}{2}}), \exp(-n^{\frac{\delta}{2}})]$ is
also large enough, which is a contradiction. Therefore, there exists
$j^-\in\Z$ such that
\begin{align*} 0\leq
x_{j^--1}<x_{j^-}<x_{j^-+1}<\exp(-n^{\frac{\delta}{2}}).\end{align*}

Similarly, there exists $j^+\in\Z$ such that
\begin{align*} 1-\exp(-n^{\frac{\delta}{2}})\leq
x_{j^+-1}<x_{j^+}<x_{j^++1}<1.\end{align*}The proof of Lemma
\ref{lbspace} is completed.\End

\begin{Remark}
From the proof of Lemma \ref{lbspace}, it is easy to see that each
of $[0, \exp(-n^{\frac{\delta}{2}})]$ and
$[1-\exp(-n^{\frac{\delta}{2}}), 1]$ contains a large number of
points of the minimal configuration $(x_i)_{i\in\Z}$ for $n$ large
enough.
\end{Remark}

 By Lemma \ref{count w} and Lemma \ref{lbspace}, without loss of
generality, one can assume that
\begin{equation}\label{jj}j^+-j^-\geq
C\left(n^{\frac{a}{2}+\frac{2\delta}{3}}\right).\end{equation}

 If $\xi\in
\mathcal {A}_\omega^{h_n}$, then $P_\omega^{h_n}(\xi)=0$. Hence, it
suffices to consider the case with $\xi\not\in \mathcal
{A}_\omega^{h_n}$ for destruction of invariant circles.  Let
$(\xi^-,\xi^+)$ be the complementary interval of $\mathcal
{A}_\omega^{h_n}$ in $\R$ and contains $\xi$. Let
$\pmb{\xi^{\pm}}=(\xi_i^{\pm})_{i\in\Z}$ be the minimal
configurations with rotation symbol $\omega$ satisfying
$\xi_0^{\pm}=\xi^{\pm}$ and  let $(\xi_i)_{i\in \Z}$ be a minimal
configuration of $h_n$ with rotation symbol $\omega$ satisfying
$\xi_0=\xi$ and $\xi_i^-\leq \xi_i\leq \xi_i^+$. By the definition
of Peierls barrier, we have
\[P_\omega^{h_n}(\xi)=\sum_{i\in \Z}(h_n(\xi_i,\xi_{i+1})-h_n(\xi_i^-,\xi_{i+1}^-)).\]
Since $P_\omega^{h_n}(\xi)$ is 1-periodic with respect to $\xi$,
without loss of generality, we assume that $\xi\in [0,1]$. We set
$d(x)=\min\{|x|,|x-1|\}$ and denote $\exp(-n^{\frac{\delta}{2}})$ by
$\epsilon(n)$. By Lemma \ref{lbspace}, there exist $i^-,\ i^+$ such
that
\begin{equation}\label{5}
d(\xi_i^-)<\epsilon(n)\quad \text{and}\quad
\xi_{i+1}^--\xi_{i-1}^-\leq\epsilon(n)\quad \text{for}\quad i=i^-,\
i^+. \end{equation}
 Thanks to Aubry's crossing lemma, we have $\xi_i^-\leq
\xi_i\leq \xi_i^+\leq \xi_{i+1}^-$. Hence,
\[\xi_i-\xi_i^-\leq\epsilon(n)\quad \text{for}\quad i=i^-,\ i^+.\]

\subparagraph{Proof of lemma \ref{apprx}} In the following, we will
prove Lemma \ref{apprx} with the method similar to the one developed
by Mather in [M3]. The proof can be proceeded in the following two
steps.

\noindent \textbf{Step 1} We consider the number of the elements in
a segment of the configuration as the length of the segment. In the
first step,
 we approximate
$P_\omega^{h_n}(\xi)$ for $\xi\in
\left[\frac{1}{2}-\frac{1}{n^a},\frac{1}{2}+\frac{1}{n^a}\right]$ by
the difference of the actions of the segments of length $i^+-i^-+1$.
To achieve that, we define the following configurations
\begin{align*}
x_i=\left\{\begin{array}{ll}\hspace{-0.4em} \xi_i,& i\neq i^-,\ i^+,\\
\hspace{-0.4em}\xi_i^-,& i=i^-,\ i^+,
\end{array}\right.
\quad \text{and}\ \quad y_i=\left\{\begin{array}{ll}
\hspace{-0.4em}\xi_i,& i^-<i<i^+,\\
\hspace{-0.4em}\xi_i^-,& i\leq i^-,\ i\geq i^+,
\end{array}\right.
\end{align*}where $(\xi_i)_{i\in\Z}$ is a minimal configuration.
It is easy to see that $\xi_0=\xi$ is contained both of
$(x_i)_{i\in\Z}$ and $(y_i)_{i\in\Z}$ up to the rearrangement of the
index $i$ since $\xi\in
\left[\frac{1}{2}-\frac{1}{n^a},\frac{1}{2}+\frac{1}{n^a}\right]$.
Hence, by the minimality of $(\xi_i)_{i\in \Z}$ satisfying
$\xi_0=\xi$, we have
\begin{equation}\label{c}
P_\omega^{h_n}(\xi)\leq\sum_{i\in
\Z}(h_n(y_i,y_{i+1})-h_n(\xi_i^-,\xi_{i+1}^-)). \end{equation} Since
$\omega$ is irrational, then $(x_i)_{i\in \Z}$ is asymptotic to
$(\xi_i^-)_{i\in \Z}$, which together with the minimality of
$(\xi_i^-)_{i\in \Z}$ yields
\begin{equation}\label{b}
\sum_{i\in
\Z}(h_n(y_i,y_{i+1})-h_n(\xi_i^-,\xi_{i+1}^-)\leq\sum_{i\in
\Z}(h_n(x_i,x_{i+1})-h_n(\xi_i^-,\xi_{i+1}^-)). \end{equation} We
set
\[h(x_i,...,x_j)=\sum_{i\leq s<j}h(x_s,x_{s+1}),\] then
\[\sum_{i\in
\Z}(h_n(x_i,x_{i+1})-h_n(\xi_i,\xi_{i+1}))=\sum_{i=i^-,i^+}(h_n(\xi_{i-1},\xi_i^-,\xi_{i+1})-h_n(\xi_{i-1},\xi_i,\xi_{i+1})).\]
 By the construction of $v_n$ and Lemma \ref{lbspace}, we have $v_n(\xi_{i^-}),\ v_n(\xi_{i^-}^-)=0$. It
follows that
\begin{align*}
h_n(&\xi_{i^--1},\xi_{i^-}^-, \xi_{i^-+1})-h_n(\xi_{i^--1},\xi_{i^-},\xi_{i^-+1})\\
&=h_n(\xi_{i^--1},\xi_{i^-}^-)+h_n(\xi_{i^-}^-,\xi_{i^-+1})-h_n(\xi_{i^--1},\xi_{i^-})-h_n(\xi_{i^-},\xi_{i^-+1}),\\
&=(\xi_{i^-}-\xi_{i^-}^-)(\xi_{i^--1}+\xi_{i^-}^-+\xi_{i^-}+\xi_{i^-+1})+u_n(\xi_{i^-})-u_n(\xi_{i^-}^-),\\
&\leq 4(\xi_{i^-}-\xi_{i^-}^-)\epsilon(n)+u_n'(\eta)(\xi_{i^-}-\xi_{i^-}^-),\\
&\leq 4\epsilon(n)^2+\frac{2\pi}{n^a}\sin(2\pi\eta)\epsilon(n),\\
&\leq C\epsilon(n)^2,
\end{align*}
where $\eta\in (\xi_{i^-},\xi_{i^-}^-)$. It is similar to obtain
\[h_n(\xi_{i^+-1},\xi_{i^+}^-,\xi_{i^++1})-h_n(\xi_{i^+-1},\xi_{i^+},\xi_{i^++1})\leq C\epsilon(n)^2.\]
Hence,
\begin{equation}\label{a}
\sum_{i\in \Z}(h_n(x_i,x_{i+1})-h_n(\xi_i,\xi_{i+1}))\leq
C\epsilon(n)^2.\end{equation} Moreover,
\begin{align*}
\sum_{i\in \Z}(h_n&(x_i,x_{i+1})-h_n(\xi_i^-,\xi_{i+1}^-))\\
&=\sum_{i\in
\Z}(h_n(x_i,x_{i+1})-h_n(\xi_i,\xi_{i+1})+h_n(\xi_i,\xi_{i+1})-h_n(\xi_i^-,\xi_{i+1}^-)),\\
&=\sum_{i\in
\Z}(h_n(x_i,x_{i+1})-h_n(\xi_i,\xi_{i+1}))+P_\omega^{h_n}(\xi),\\
&\leq P_\omega^{h_n}(\xi)+C\epsilon(n)^2.
\end{align*}
Therefore, it follows from $(\ref{c})$ and $(\ref{b})$ that
\begin{equation}\label{step1}
P_\omega^{h_n}(\xi)\leq\sum_{i\in
\Z}(h_n(y_i,y_{i+1})-h_n(\xi_i^-,\xi_{i+1}^-))\leq
P_\omega^{h_n}(\xi)+C\epsilon(n)^2,
\end{equation}
where \begin{equation}\label{lnew}\sum_{i\in
\Z}(h_n(y_i,y_{i+1})-h_n(\xi_i^-,\xi_{i+1}^-))=h_n(y_{i^-},...,y_{i^+})-h_n(\xi_{i^-}^-,...,\xi_{i^+}^-).\end{equation}

\noindent \textbf{Step 2} It follows from [M4] that the Peierls's
barrier $P_{0^+}^{h_n}(\xi)$ could be defined as follows
\begin{equation}\label{t}
P_{0^+}^{h_n}(\xi)=\min_{\eta_0=\xi}\sum_{i\in
\Z}h_n(\eta_i,\eta_{i+1})-\min\sum_{i\in \Z}h_n(z_i,z_{i+1}),
\end{equation}
where $(\eta_i)_{i\in\Z}$ and  $(z_i)_{i\in\Z}$ are monotone
increasing configurations limiting on $0,\ 1$. We set
\begin{equation*}
\begin{cases}
K(\xi)=\min_{\eta_0=\xi}\sum_{i\in \Z}h_n(\eta_i,\eta_{i+1}),\\
K=\min\sum_{i\in \Z}h_n(z_i,z_{i+1}).
\end{cases}
\end{equation*}

First of all, it is easy to see that $K(\xi)$ and $K$ are bounded.
Second, $P_{0^+}^{h_n}(\xi)=0$ for $\xi=0$ or $1$. Hence, we only
need to consider the case with $\xi\in (0,1)$. Following the ideas
of [M6], let $\pmb{\xi^-}$ and $\pmb{\xi^+}$ be minimal
configurations of rotation symbol $0^+$ and let $(\xi_0^-,\xi_0^+)$
be the complementary interval of $\mathcal {A}_{0^+}^{h_n}$ and
contains $\xi$. Based on the definition
\[P_{0^+}^{h_n}(\xi)=\min_{x_0=\xi}\{G_{0^+}(\bold{x})|\xi_i^-\leq \zeta_i\leq\xi_i^+\},\]
where \[G_{0^+}(\pmb{\zeta})=\sum_{i\in
\Z}(h_n(\zeta_i,\zeta_{i+1})-h_n(\xi_i^-,\xi_{i+1}^-))=-K+\sum_{i\in
\Z}h_n(\zeta_i,\zeta_{i+1}),\] the proof of $(\ref{t})$ will be
completed when we verify that the configuration $(\zeta_i)_{i\in
\Z}$ achieving the minimum in the definition of $K(\xi)$ satisfies
$\xi_i^-\leq \zeta_i\leq\xi_i^+$. It can be easily obtained by
Aubry's crossing lemma. In fact, since $(\xi_i^-)_{i\in\Z}$ and
$(\zeta_i)_{i\in\Z}$ are minimal and both are $\alpha$-asymptotic to
$0$ as well as $\omega$-asymptotic to $1$, by Aubry's crossing
lemma, $(\xi_i^-)_{i\in\Z}$ and $(\zeta_i)_{i\in\Z}$ do not cross.
Similarly $(\xi_i^+)_{i\in\Z}$ and $(\zeta_i)_{i\in\Z}$ do not
cross. It follows from $\zeta_0\in (\xi_0^-,\xi_0^+)$ that
$(\zeta_i)_{i\in \Z}$ achieving the minimum in the definition of
$K(\xi)$ satisfies $\xi_i^-\leq \zeta_i\leq\xi_i^+$.

In the following, we will compare $K$, $K(\xi)$ with
$h_n(\xi_{i^-}^-,...,\xi_{i^+}^-)$, $h_n(y_{i^-},...,y_{i^+})$
respectively, here $h_n(\xi_{i^-}^-,...,\xi_{i^+}^-)$ and
$h_n(y_{i^-},...,y_{i^+})$ are as the same as the notations in
(\ref{lnew}).

First, we consider $K$ and $h_n(\xi_{i^-}^-,...,\xi_{i^+}^-)$. Let
$(z_i)_{i\in\Z}$ be a monotone increasing configuration limiting on
$0,\ 1$ such that $K=\sum_{i\in \Z}h_n(z_i,z_{i+1})$. By Lemma
\ref{count 0},
\[\sharp\{i\in
\Z|z_i\in [\epsilon(n),1-\epsilon(n)]\}\leq
Cn^{\frac{a}{2}+\frac{\delta}{2}}.\] On the other hand, since
$(z_i)_{i\in\Z}$ has the rotation number $0^+$, then from
(\ref{jj}), it follows that up to the rearrangement of the index
$i$, there exists a subset of length $i^+-i^-$ of $(z_i)_{i\in\Z}$,
denoted by $\{z_{i^-},z_{i^-+1},\ldots,z_{i^+-1},z_{i^+}\}$ such
that
\[z_{i^-+1}\leq\epsilon(n),\quad z_{i^+-1}\geq 1-\epsilon(n).\]

By Lemma \ref{lbspace}, $\xi_{i^-+1}^->0$ and $\xi_{i^+-1}^-<1$ for
the minimal configuration $(\xi_i^-)_{i\in\Z}$. We consider the
configuration $(\bar{x}_i)_{i\in\Z}$ defined by
\begin{equation*}
\left\{\begin{array}{ll}
\hspace{-0.4em}\bar{x}_i=\xi_i^-,&i^-<i<i^+,\\
\hspace{-0.4em}\bar{x}_i=0,& i\leq i^-,\\
\hspace{-0.4em}\bar{x}_i=1,& i\geq i^+.
\end{array}\right.
\end{equation*}
By the definition of $h_n$, $h_n(\bar{x}_i,\bar{x}_{i+1})=0$ for
$i<i^-$ or $i \geq i^+$, then
\[\sum_{i\in\Z}h_n(\bar{x}_i,\bar{x}_{i+1})=h_n(\bar{x}_{i^-},...,\bar{x}_{i^+}).\]
Moreover, by the minimality of $(z_i)_{i\in\Z}$, we have
\begin{equation}\label{aa}
K\leq
\sum_{i\in\Z}h_n(\bar{x}_i,\bar{x}_{i+1})=h_n(\bar{x}_{i^-},...,\bar{x}_{i^+}).
\end{equation}
By the construction of $h_n$, we have
$v_n(\bar{x}_{i^-+1})=v_n(\xi_{i^-+1}^-)=0$. Hence,
\begin{equation}\label{55}
\begin{split}
h_n(&\bar{x}_{i^-},\bar{x}_{i^-+1})-h_n(\xi_{i^-}^-,\xi_{i^-+1}^-)\\
&=\frac{1}{2}(\bar{x}_{i^-}-\bar{x}_{i^-+1})^2+u_n(\bar{x}_{i^-+1})
-\frac{1}{2}(\xi_{i^-}^--\xi_{i^-+1}^-)^2-u_n(\xi_{i^-+1}^-),\\
&=\frac{1}{2}(\xi_{i^-+1}^-)^2-\frac{1}{2}(\xi_{i^-+1}^--\xi_{i^-}^-)^2,\\
&=\frac{1}{2}\xi_{i^-}^-(2\xi_{i^-+1}^--\xi_{i^-}^-),\\
&\ \leq C\epsilon(n)^2.
\end{split}
\end{equation}
It is similar to obtain
\begin{equation}\label{3}
h_n(\bar{x}_{i^+-1},\bar{x}_{i^+})
-h_n(\xi_{i^+-1}^-,\xi_{i^+}^-)\leq C\epsilon(n)^2. \end{equation}
 Since
\begin{align*}
h_n(&\bar{x}_{i^-},...,\bar{x}_{i^+})-h_n(\xi_{i^-}^-,...,\xi_{i^+}^-)\\
&=h_n(\bar{x}_{i^-},\bar{x}_{i^-+1})-h_n(\xi_{i^-}^-,\xi_{i^-+1}^-)+h_n(\bar{x}_{i^+-1},\bar{x}_{i^+})
-h_n(\xi_{i^+-1}^-,\xi_{i^+}^-),
\end{align*}then
\begin{equation}\label{aaa}
h_n(\bar{x}_{i^-},...,\bar{x}_{i^+})-h_n(\xi_{i^-}^-,...,\xi_{i^+}^-)\leq
C\epsilon(n)^2.
\end{equation}
 From $(\ref{aa})$ and $(\ref{aaa})$ we have
\begin{equation}\label{1}
 K\leq h_n(\xi_{i^-}^-,...,\xi_{i^+}^-)+C\epsilon(n)^2.
\end{equation}

To obtain the reverse inequality of $(\ref{1})$, we consider the
configuration as follows
\begin{equation*}
\left\{\begin{array}{ll}
\hspace{-0.4em}\tilde{x}_i=z_i,&i^-<i<i^+,\\
\hspace{-0.4em}\tilde{x}_i=0,&i\leq i^-,\\
\hspace{-0.4em}\tilde{x}_i=1,&i\geq i^+.
\end{array}\right.
\end{equation*}

From the definition of $h_n$, it follows that $v_n(z_{i^-+1})=0$ and
$h_n(z_i,z_{i+1})\geq 0$ for all $i\in \Z$. Moreover, we have
\begin{align*}
h_n(&\tilde{x}_{i^-},...,\tilde{x}_{i^+})-K\\
&=h_n(\tilde{x}_{i^-},\tilde{x}_{i^-+1})+h_n(\tilde{x}_{i^+-1},\tilde{x}_{i^+})-\sum_{i<i^-,i\geq
i^+}h_n(z_i,z_{i+1}),\\
&\leq\frac{1}{2}(z_{i^-+1})^2+u_n(z_{i^-+1})+\frac{1}{2}(z_{i^+-1}-1)^2,\\
&\leq u_n'(\eta)z_{i^-+1}+C_1\epsilon(n)^2,\\
&\leq 2\pi n^{-a}\sin(2\pi\eta)z_{i^-+1}+C_1\epsilon(n)^2,\\
&\leq C_2n^{-a}(z_{i^-+1})^2+C_1\epsilon(n)^2,\\
&\leq C\epsilon(n)^2.
\end{align*}
where $\eta\in (0,z_{i^-+1})$. Namely
\begin{equation}\label{bb}
h_n(\tilde{x}_{i^-},...,\tilde{x}_{i^+})\leq K+C\epsilon(n)^2.
\end{equation}
Furthermore, we consider the finite segment of the configuration
defined by
\begin{equation*}
\left\{\begin{array}{ll}
\hspace{-0.4em}\eta_i=\tilde{x}_i,&i^-<i<i^+,\\
\hspace{-0.4em}\eta_i=\xi_i^-,&i=i^-,\\
\hspace{-0.4em}\eta_i=\xi_i^-,& i=i^+.
\end{array}\right.
\end{equation*}
Then, the minimality of $(\xi_i^-)_{i\in\Z}$ implies
$h_n(\xi_{i^-}^-,...,\xi_{i^+}^-)\leq
h_n(\eta_{i^-},...,\eta_{i^+})$. Hence, by $(\ref{bb})$, we have
\begin{equation}\label{cc}
h_n(\xi_{i^-}^-,...,\xi_{i^+}^-)\leq
K+C\epsilon(n)^2+h_n(\eta_{i^-},...,\eta_{i^+})-h_n(\tilde{x}_{i^-},...,\tilde{x}_{i^+}),
\end{equation}
where
\begin{align*}
 h_n(&\eta_{i^-},...,\eta_{i^+})-h_n(\tilde{x}_{i^-},...,\tilde{x}_{i^+})\\
 &=h_n(\eta_{i^-},\eta_{i^-+1})-h_n(\tilde{x}_{i^-},\tilde{x}_{i^-+1})+h_n(\eta_{i^+-1},\eta_{i^+})-h_n(\tilde{x}_{i^+-1},\tilde{x}_{i^+}).
\end{align*}
By the deduction as similar as $(\ref{55})$, we have
\begin{align*}
 &h_n(\eta_{i^-},\eta_{i^-+1})-h_n(\tilde{x}_{i^-},\tilde{x}_{i^-+1})\leq C\epsilon(n)^2,\\
 &h_n(\eta_{i^+-1},\eta_{i^+})-h_n(\tilde{x}_{i^+-1},\tilde{x}_{i^+})\leq C\epsilon(n)^2.
\end{align*}
Moreover, \begin{equation}\label{dd}
 h_n(\eta_{i^-},...,\eta_{i^+})-h_n(\tilde{x}_{i^-},...,\tilde{x}_{i^+})\leq
 C\epsilon(n)^2.
\end{equation}
 Hence, from $(\ref{cc})$ and $(\ref{dd})$, it follows that
\begin{equation}\label{2}
h_n(\xi_{i^-}^-,...,\xi_{i^+}^-)\leq K+C\epsilon(n)^2,
\end{equation}
which together with $(\ref{1})$ implies
\begin{equation}\label{13}
|h_n(\xi_{i^-}^-,...,\xi_{i^+}^-)-K|\leq C\epsilon(n)^2.
\end{equation}

Next, we compare $h_n(y_{i^-},...,y_{i^+})$ with $K(\xi)$.

Since $(\xi_i)_{i\in\Z}$ is minimal among all configurations with
rotation symbol $\omega$ satisfying $\xi_0=\xi$. By $(\ref{5})$ and
Aubry's crossing lemma, we have \[d(\xi_i)\leq \epsilon(n),\quad
\text{for}\ i=i^-,i^+,\] where $d(\xi_i)=\min\{|\xi_i|,|\xi_i-1|\}$.
By an argument as similar as the one in the comparison between $K$
and $h_n(\xi_{i^-}^-,...,\xi_{i^+}^-)$, we have
\begin{equation}\label{12}
|h_n(\xi_{i^-},...,\xi_{i^+})-K(\xi)|\leq C\epsilon(n)^2.
\end{equation}
By the construction of $(y_i)_{i\in\Z}$, namely
\begin{equation*}
y_i=\left\{\begin{array}{ll}
\hspace{-0.4em}\xi_i,& i^-<i<i^+,\\
\hspace{-0.4em}\xi_i^-,& i\leq i^-,\ i\geq i^+,
\end{array}\right.
\end{equation*}
we have
\begin{align*}
h_n(&y_{i^-},...,y_{i^+})-h_n(\xi_{i^-},...,\xi_{i^+})\\
&=h_n(\xi_{i^-}^-,\xi_{i^-+1})-h_n(\xi_{i^-},\xi_{i^-+1})+h_n(\xi_{i^+-1},\xi_{i^+}^-)-h_n(\xi_{i^+-1},\xi_{i^+}).
\end{align*}
By the deduction as similar as $(\ref{aaa})$, we have
\begin{equation}\label{11}
|h_n(y_{i^-},...,y_{i^+})-h_n(\xi_{i^-},...,\xi_{i^+})|\leq
C\epsilon(n)^2.
\end{equation}

Finally, from $(\ref{step1})$, $(\ref{13})$, $(\ref{12})$ and
$(\ref{11})$ we obtain
\begin{align*}
|P_{\omega}^{h_n}(\xi)-P_{0^+}^{h_n}(\xi)|&\leq|h_n(y_{i^-},...,y_{i^+})-h_n(\xi_{i^-}^-,...,\xi_{i^+}^-)+K-K(\xi)|+C_1\epsilon(n)^2,\\
&\leq
|h_n(\xi_{i^-},...,\xi_{i^+})-K(\xi)|+|h_n(\xi_{i^-}^-,...,\xi_{i^+}^-)-K|\\
&\ +|h_n(y_{i^-},...,y_{i^+})-h_n(\xi_{i^-},...,\xi_{i^+})|+C_1\epsilon(n)^2,\\
&\leq C\epsilon(n)^2,\\
&=C\exp\left(-2n^{\frac{\delta}{2}}\right),
\end{align*}
which completes the proof of Lemma \ref{apprx}.\End

\paragraph{Proof of Theorem \ref{MR}} Based on the preparation above, it
is easy to prove Theorem \ref{MR}. We assume that there exists an
invariant circle with rotation number
$0<\omega<n^{-\frac{a}{2}-\delta}$ for $h_n$, then
$P_\omega^{h_n}(\xi)\equiv 0$ for every $\xi\in \R$. By Lemma
\ref{apprx}, we have
\begin{equation}\label{pp}
|P_{0^+}^{h_n}(\xi)|\leq
C\exp\left(-2n^{\frac{\delta}{2}}\right),\quad\text{for}\quad \xi\in
\left[\frac{1}{2}-\frac{1}{n^a},\frac{1}{2}+\frac{1}{n^a}\right].
\end{equation}

On the other hand, $(\ref{lowbound})$ implies that there exists a
point $\tilde{\xi}\in
\left[\frac{1}{2}-\frac{1}{n^a},\frac{1}{2}+\frac{1}{n^a}\right]$
such that
\[P_{0^+}^{h_n}(\tilde{\xi})\geq n^{-s}.\]Hence, we have
\[n^{-s}\leq C\exp\left(-2n^{\frac{\delta}{2}}\right).\] It is an obvious
contradiction for $n$ large enough. Therefore, there exists no
invariant circle with rotation number
$0<\omega<n^{-\frac{a}{2}-\delta}$.

For $-n^{-\frac{a}{2}-\delta}<\omega<0$, by comparing
$P_{\omega}^{h_n}(\xi)$ with $P_{0^-}^{h_n}(\xi)$, the proof is
similar. We omit the details. Therefore, the proof of Theorem
\ref{MR} is completed. \End
\subsubsection{$C^\omega$ case}
We will prove the following theorem:
\begin{Theorem}\label{maint1} Given an integrable generating function $h_0$, a rotation number $\omega$ and a small positive constant
$\delta$, there exists a sequence of real-analytic $(h_n)_{n\in\N}$
such that $h_n\rightarrow h_0$ in the $C^{3-\delta}$ topology  and
the exact monotone area-preserving twist maps generated by
$(h_n)_{n\in\N}$ admit no invariant circles with the rotation number
$\omega$. \end{Theorem}
\paragraph{Construction of the generating functions}
 Consider a completely integrable system with
the generating function
\[h_0(x,x')=\frac{1}{2}(x-x')^2, \quad x,x'\in \R.\] We construct the perturbation consisting of two
parts. The first one is
\begin{equation}\label{31}
u_n(x)=\frac{1}{n^a}(1-\cos(2\pi x) ),\quad x\in \R,\end{equation}
where $n\in \N$ and $a$ is a positive constant independent of $n$.

We construct  the second part of the perturbation in the following.
Let $p_N(x)$ be a trigonometric polynomial of degree $N$. It is easy
to see that for any $r>0$,
\begin{equation}
||p_N(x)||_r\leq e^{rN}||p_N(x)||,
\end{equation}
where $||p_N(x)||_r$ denotes the maximum of $|p_N(z)|$ in the strip
$S_r=\{z\in\mathbb{C}|\,|\text{Im}z|\leq r\}$ of width $2r$ in the
complex plane and $||p_N(x)||$ denotes the maximum of $|p_N(x)|$ on
the real line. Without loss of generality, we take $r=1$, a.e.
\begin{equation}
||p_N(x)||_1\leq e^{N}\max|p_N(x)|.
\end{equation}
Then, by Cauchy estimates, for any fixed $s>0$, we have
\begin{equation}\label{cachy}
||p_N(x)||_{C^s}\leq C_se^N\max|p_N(x)|,
\end{equation}
where $C_s$ is a constant depending only on $s$.

Based on Lemma \ref{lowstep}, we need to construct a real analytic
function with a ``bump" in correspondence with the interval
$\Lambda_n$ satisfying
\begin{equation}\label{length}
\mathcal {L}(\Lambda_n)\sim n^{-\frac{a}{2}}\quad\text{and}\quad
\Lambda_n\subset\left[\frac{1}{4},\frac{3}{4}\right],
\end{equation}
where $\mathcal {L}(\Lambda_n)$ denotes the Lebesgue measure of
$\Lambda_n$ and $f\sim g$ means that $\frac{1}{C}g<f<C g$ holds for
some constant $C>1$.

The ``bump" will be accomplished by using Jackson's approximation
theorem (see [Z, p115]). It states that let $\phi(x)$ be an
$k$-times differentiable periodic function on $\R$, then for every
$N\in\N$, there exists a trigonometric polynomial $p_N(x)$ of degree
$N$ such that
\[\max|p_N(x)-\phi(x)|\leq A_kN^{-k}||\phi(x)||_{C^k},\]
where $A_k$ is a constant depending only on $k\in\N$.

We take a $C^\infty$ bump function $\phi$ supported on the interval
$\Lambda_n$, whose maximum is equal to $2$. By (\ref{length}), the
length of $\Lambda_n$ is bounded by $Cn^{-\frac{a}{2}}$, one can
choose $\phi(x)$ such that
\begin{equation}
||\phi(x)||_{C^k}\sim
\left(n^{\frac{a}{2}}\right)^k=n^{\frac{ak}{2}}.
\end{equation}
Then, choosing $N$ large enough to achieve
\begin{equation}\label{N}
A_k N^{-k}||\phi(x)||_{C^k}<\sigma\ll 1,
\end{equation}
where $\sigma$ will be determined in the following, by Jackson's
approximation theorem, we can construct a trigonometric polynomial
$p_N(x)$ of degree $N$ such that:
\begin{equation}
\left\{\begin{array}{ll}\hspace{-0.4em}\max p_N(x)\geq 1,&\text{attained on}\ \Lambda_n,\\
\hspace{-0.4em}|p_N(x)|\leq\sigma,&\text{on}\ [0,1]\backslash\Lambda_n.\\
\end{array}\right.
\end{equation}
By (\ref{N}), we have
\begin{equation}\label{nn}
N\geq C\sigma^{-\frac{1}{k}}n^{\frac{a}{2}}.
\end{equation}
Finally, we consider the normalized trigonometric polynomial
\begin{equation}
\tilde{p}_N(x)=e^{-2N}\left(\frac{p_N(x)}{\max p_N(x)}\right)^2.
\end{equation}
From (\ref{cachy}), $\tilde{p}_N(x)$ satisfies:
\begin{equation}
\left\{\begin{array}{ll}\hspace{-0.4em}\tilde{p}_N(x)\geq 0,\\
\hspace{-0.4em}||\tilde{p}_N(x)||_{C^s}\leq C,\\
\hspace{-0.4em}\max \tilde{p}_N(x)=e^{-2N},&\text{attained on}\ \Lambda_n,\\
\hspace{-0.4em}|\tilde{p}_N(x)|\leq \sigma^2e^{-2N},&\text{on}\ [0,1]\backslash\Lambda_n.\\
\end{array}\right.
\end{equation}
Based on preparations above, we can construct the second part of the
perturbation as follow
\begin{equation}
v_n(x)=u_n(x)\tilde{p}_N(x)=\frac{1}{n^a}(1-\cos 2\pi
x)\tilde{p}_N(x).
\end{equation}
It is easy to see $v_n$ satisfies the following properties:
\begin{equation}\label{vn}
\left\{\begin{array}
{ll}\hspace{-0.4em}v_n(x)\geq 0,&\\
\hspace{-0.4em}||v_n(x)||_{C^s}\leq Cn^{-a},&\\
\hspace{-0.4em}\max v_n(x)\geq e^{-2N}n^{-a},&\text{attained on}\ \Lambda_n,\\
\hspace{-0.4em}|v_n(x)|\leq C\sigma^2e^{-2N}n^{-a},&\text{on}\
[0,1]\backslash\Lambda_n.
\end{array}\right.
\end{equation}

So far, we complete the construction of the generating function of
the nearly integrable system, a.e.
\begin{equation}\label{h}
h_n(x,x')=h_0(x,x')+u_n(x')+v_n(x'),
\end{equation}
where $n\in\N$.

\paragraph{Proof of Theorem \ref{maint1}}
 If $\omega\in \Q$, then the invariant circles with rotation
number $\omega$ could be easily destructed by an analytic
perturbation arbitrarily close to $0$. Therefore it suffices to
consider the irrational $\omega$. Firstly, we prove the
non-existence of invariant circles with a small enough rotation
number. More precisely, we have the following Lemma:
\begin{Lemma}\label{MRr} For $\omega\in\R\backslash\Q$ and $n$ large enough,
the exact area-preserving monotone twist map generated by $h_n$
admits no invariant circle with the rotation number satisfying
\[|\omega|<n^{-a-\delta}, \] where $\delta$ is a small positive constant independent of $n$.
\end{Lemma}
\Proof First of all, we estimate the lower bound of $P_{0^+}^{h_n}$.
 Let $(\xi_i)_{i\in \Z}$
be a minimal configuration of $h_n$ defined by (\ref{h}) with
rotation symbol $0^+$ satisfying $\xi_0=\eta$, where $\eta$
satisfies $v_n(\eta)=\max v_n(x)$ and let $(x_i)_{i\in \Z}$ be the
minimal configuration of
$\bar{h}_n(x_i,x_{i+1})=h_0(x_i,x_{i+1})+u_n(x_{i+1})$ with rotation
symbol $0^+$ satisfying $x_0\in [0,1]\backslash\Lambda_n$, then
\begin{align*}
\sum_{i\in
\Z}(h_n(&\xi_i,\xi_{i+1})-h_n(\xi_i^-,\xi_{i+1}^-))\\
&\geq v_n(\eta)+\sum_{i\in \Z}\bar{h}_n(\xi_i,\xi_{i+1})-\sum_{i\in
\Z}h_n(\xi_i^-,\xi_{i+1}^-),\\
&\geq v_n(\eta)+\sum_{i\in \Z}\bar{h}_n(x_i,x_{i+1})-\sum_{i\in
\Z}h_n(x_i,x_{i+1}),\\
&=v_n(\eta)-\sum_{i\in \Z}v_n(x_{i+1}).\\
\end{align*}
Therefore, we have shown:
\begin{equation}
P_{0^+}^{h_n}(\eta)=\min_{\xi_0=\eta}\sum_{i\in
\Z}(h_n(\xi_i,\xi_{i+1})-h_n(\xi_i^-,\xi_{i+1}^-))\geq
v_n(\eta)-\sum_{i\in \Z}v_n(x_{i+1}).
\end{equation}
By (\ref{vn}), we have
\begin{equation}
v_n(\eta)\geq e^{-2N}n^{-a}.
\end{equation}
It follows that
\begin{equation}
\sum_{i\in \Z}v_n(x_{i+1})\leq
\sigma^2e^{-2N}\sum_{i\in\Z}u_n(x_{i+1})\leq\sigma^2e^{-2N}\sum_{i\in\Z}\frac{1}{4}(x_{i+1}-x_{i-1})^2\leq\sigma^2e^{-2N}.
\end{equation}
Hence,
\begin{equation}
P_{0^+}^{h_n}(\eta)\geq e^{-2N}(n^{-a}-\sigma^2),
\end{equation}
we choose then $\sigma$ (consequently $N$) in such a way that
\[\frac{1}{4}n^{-a}-\sigma^2\geq 0,\]
which implies
\[\sigma\leq\frac{1}{2}n^{-\frac{a}{2}}.\]
By (\ref{nn}), if follows that
\begin{equation}
N\geq Cn^{\frac{a}{2}+\frac{a}{2k}}.
\end{equation}
Therefore,
\begin{equation}\label{lowb}
\max_NP_{0^+}^{h_n}(\eta)\geq
n^{-a}\exp\left(-Cn^{\frac{a}{2}+\frac{a}{2k}}\right).
\end{equation}

 Second, following a similar argument as \cite{W1}, we have
\begin{equation}\label{app}
|P_{\omega}^{h_n}(\xi)-P_{0^+}^{h_n}(\xi)|\leq
C\exp\left(-2n^{\frac{a}{2}+\frac{\delta}{2}}\right).\end{equation}where
$\xi\in \Lambda_n$ and $\delta$ is a small positive constant
independent of $n$. Here $\Lambda_n$ is as the same as the notation
in (\ref{length}).

Based on the preparations above, it is easy to prove Lemma
\ref{MRr}. We assume that there exists an invariant circle with
rotation number $0<\omega<n^{-a-\delta}$ for $h_n$, then
$P_\omega^{h_n}(\xi)\equiv 0$ for every $\xi\in \R$. By (\ref{app}),
we have
\begin{equation}\label{pp}
|P_{0^+}^{h_n}(\xi)|\leq
C\exp\left(-2n^{\frac{a}{2}+\frac{\delta}{2}}\right),\quad\text{for}\quad
\xi\in \Lambda_n.
\end{equation}

On the other hand, $(\ref{lowb})$ implies that there exists a point
$\eta\in \Lambda_n$ such that
\[P_{0^+}^{h_n}(\eta)\geq n^{-a}\exp\left(-Cn^{\frac{a}{2}+\frac{a}{2k}}\right).\]Hence, we have
\begin{equation}\label{cotr}
n^{-a}\exp\left(-Cn^{\frac{a}{2}+\frac{a}{2k}}\right)\leq
C\exp\left(-2n^{\frac{a}{2}+\frac{\delta}{2}}\right).
\end{equation}
To achieve the contradiction, it suffices to take
\[k>\frac{3a}{2\delta},\]
which implies
\[\frac{a}{2k}<\frac{\delta}{3}<\frac{\delta}{2}.\]
Hence, for $n$ large enough
\begin{equation}
n^{-a}\exp\left(-Cn^{\frac{a}{2}+\frac{a}{2k}}\right)\geq
C\exp\left(-2n^{\frac{a}{2}+\frac{\delta}{2}}\right),
\end{equation}
which contradicts (\ref{cotr}).
 Therefore,
there exists no invariant circle with rotation number
$0<\omega<n^{-a-\delta}$.

For $-n^{-a-\delta}<\omega<0$, by comparing $P_{\omega}^{h_n}(\xi)$
with $P_{0^-}^{h_n}(\xi)$, the proof is similar. We omit the
details. This completes the proof of Lemma \ref{MRr}. \End

By Lemma \ref{Herm}, the case with a given irrational rotation
number can be easily reduced to the one with a small enough rotation
number. For the sake of simplicity of notations, we denote $Q_{q_n}$
by $Q_n$ and the same to $u_{q_n}, v_{q_n}$ and $h_{q_n}$. Let
\[Q_n(x)={q_n}^{-2}(u_n(q_nx)+v_n(q_nx)),\] where $(q_n)_{n\in \N}$
is a sequence satisfying Dirichlet approximation
\begin{equation}\label{diri}
|q_n\omega-p_n|<\frac{1}{q_n}, \end{equation} where $p_n\in \Z$ and
$q_n\in \N$. Since $\omega\in\R\backslash\Q$, we say
$q_n\rightarrow\infty$ as $n\rightarrow\infty$. Let
$\tilde{h}_n(x,x')=h_0(x,x')+Q_n(x')$, we prove Theorem \ref{maint1}
for $(\tilde{h}_n)_{n\in\N}$ as follow:

\Proof Based on Lemma \ref{MRr} and Dirichlet approximation
(\ref{diri}), it suffices to take
\[\frac{1}{q_n}\leq\frac{1}{{q_n}^{a+\delta}},\]
which implies \begin{equation}\label{33} a\leq
1-\delta.\end{equation} From the constructions of $u_n$ and $v_n$,
it follows that
\begin{align*}
 ||\tilde{h}_n&(x,x')-h_0(x,x')||_{C^r}\\
 &=||Q_n(x')||_{C^r},\\
&\leq{q_n}^{-2}(||u_n(q_nx')||_{C^r}+||v_n(q_nx')||_{C^r}),\\
&\leq{q_n}^{-2}({q_n}^{-a}(2\pi)^r{q_n}^r+C_1{q_n}^{-a}{q_n}^r),\\
&\leq C_2{q_n}^{r-a-2},
\end{align*}
where $C_1, C_2$ are positive constants only depending on $r$.

To complete the proof, it is enough to make $r-a-2<0$, which
together with (\ref{33}) implies
\[r<a+2\leq 3-\delta.\]This completes the
proof of Theorem \ref{maint1}.\End

\subsection{Case with $d\geq 3$ degrees of freedom}
\subsubsection{Preliminaries}

In $\text{T}^\ast\T^d$, a submanifold $\mathcal {T}^d$ is called
Lagrangian torus if it is diffeomorphic to the torus $\T^d$ and the
symplectic form vanishes on it. For positive definite Hamiltonian
systems, if a Lagrangian torus is invariant under the Hamiltonian
flow, it is then the graph over $\T^d$ (see \cite{BP}). An example
of Lagrangian torus is the KAM torus.
\begin{Definition}\label{dd}$\bar{\mathcal
{T}}^d$ is called $d$ dimensional KAM torus if
\begin{itemize}
\item $\bar{\mathcal {T}}^d$ is a Lipschitz graph over $\T^d$;
\item $\bar{\mathcal {T}}^d$ is invariant under the Hamiltonian flow
 $\Phi_t^H$ generated by the Hamiltonian function $H$;
\item there exists a diffeomorphism
$\phi:\ \T^d\rightarrow \bar{\mathcal {T}}^d$ such that
$\phi^{-1}\circ\Phi_H^t\circ\phi=R_\omega^t$ for any $t\in \R$,
where $R_\omega^t:\ x\rightarrow x+\omega t$ and $\omega$ is called
the rotation vector of $\bar{\mathcal {T}}^d$.
\end{itemize}
\end{Definition}

For positive definite Hamiltonian systems, each KAM torus
$\bar{\mathcal {T}}^d$ supports a minimal measure $\mu$. The
rotation number $\rho$ of $\mu$ is well defined and
$\rho(\mu)=\omega$. The rotation vector of the Lagrangian torus
 $\mathcal {T}^d$ is not well defined. If $\mathcal
{T}^d$ supports several invariant measures with different rotation
vectors. In this paper, we are only concerned with Lagrangian tori
with well defined rotation vectors.
\begin{Definition}\label{dd1}
$\mathcal {T}^d$ is called $d$ dimensional Lagrangian torus with the
rotation vector $\omega$ if
\begin{itemize}
\item $\mathcal {T}^d$ is a Lagrangian submanifold;
\item $\mathcal {T}^d$ is invariant for the Hamiltonian
flow $\Phi_H^t$ generated by $H$.
\item each orbit on $\mathcal {T}^d$ has the same rotation vector.
\end{itemize}
\end{Definition}
In \cite{H2}, it is proved that each orbit on $\mathcal {T}^d$ is an
action minimizing curve.

An arithmetic approximation of the rotation vector is found in
\cite{C}. For any given vector $\omega\in \R^d$ with $d\geq 2$,
there is a sequence of integer vectors $k_n\in\Z^d$ with
$|k_n|\rightarrow\infty$ such that
\begin{equation}\label{ap}
\left|\langle\omega, k_n\rangle\right|<\frac{C}{|k_n|^{d-1}},
\end{equation}
where $C$ is a constant independent of $n$,
\[|k|=\left(\sum_{j=1}^d k_i^2\right)^{\frac{1}{2}},\quad
\text{for}\ k=(k_1,k_2,\ldots,k_d).\]

A rotation vector $\omega\in\R^d$ is called resonant if there exists
$k\in\Z^d$ such that $\langle \omega, k\rangle=0$. Otherwise, it is
non-resonant. Obviously, a Lagrangian torus with the resonant
rotation vector can be destructed by analytic perturbation
arbitrarily close to zero. Hence, it is sufficient to consider the
Lagrangian torus with the non resonant rotation vector. In that
case, one can assume that the Lagrangian torus $\mathcal {T}^d$
supports a uniquely ergodic minimizing measure. Moreover, by
\cite{M2}, $\mathcal {T}^d$ is a Lipschitz graph over the underlying
manifold $\T^d$.
\subsubsection{$C^\infty$ case}
We will prove the following theorem:
\begin{Theorem}\label{maint} Given an integrable Hamiltonian $H_0$ with $d\ (d\geq 3)$
 degrees of freedom, a
rotation vector $\omega$ and a small positive constant $\delta$,
there exists a sequence of $C^\infty$ Hamiltonians $\{H_n\}_{n\in
\N}$ such that $H_n\rightarrow H_0$ in $C^{2d-\delta}$ topology and
 the Hamiltonian flow generated by $H_n$ does not admit the
Lagrangian torus with the rotation vector $\omega$.
\end{Theorem}

This theorem implies that the rigidity of the Lagrangian torus is as
the same as the KAM torus. Roughly speaking, the maximum of $r$ is
closely related to the arithmetic property of the rotation vector
$\omega$. If $\omega$ is a Diophantine vector,  then $r$ is at most
$2d-\delta$. If $\omega$ is a Liouville vector, then $r$ can be
 arbitrarily large. If $\omega$ can be approximated exponentially by
rational vectors, then the Lagrangian torus with the rotation vector
$\omega$ can be destructed by an arbitrarily small perturbation in
$C^\omega$ (analytic) topology (see \cite{B2}).
\paragraph{Destruction of Lagrangian torus with a special rotation
vector} The Hamiltonian function we consider here is nearly
integrable
\begin{equation}\label{ha}
H_n(q,p)=H_0(p)-P_n(q),\end{equation} where $(q,p)\in
\T^d\times\R^d$. Without loss of generality, we
assume\[H_0(p)=\frac{1}{2}|p|^2,\] for which (\ref{ha}) is a typical
mechanical system.

Since $H_n$ is strictly convex with respect to $p$, by the Legendre
transformation, the Lagrangian function corresponding to $H_n$ is
\begin{equation}\label{Lo}
L_n(q,\dot{q})=\frac{1}{2}|\dot{q}|^2+P_n(q),
\end{equation}
where $\dot{q}=\frac{\partial H_0}{\partial p}$.

Let
\[P_n(q)=\frac{1}{n^{a}}(1-\cos q_1)+v_n\left(q_1,q_2\right),\] where
$a$ is a positive constant independent of $n$. For the rotation
vector $\omega=(\omega_1,\ldots,\omega_d)$, $v_n(q_1,q_2)$ is
constructed as follow
 \begin{equation}\label{v}
\begin{cases}
 v_n\ \text{is}\ 2\pi\text{-periodic},\\
\text{supp}\,v_n\cap \{[0,2\pi]\times[-\pi,\pi]\}= B_{R_n}(q^*),\\
\max_{(q_1,q_2)\in [0,2\pi]\times[-\pi,\pi]} v_n=v_n(q_0)={|\omega_1|}^s, \\
{||v_n||}_{C^r}\sim{|\omega_1|}^{s'},
\end{cases}
\end{equation}
where $R_n=\frac{|\omega_1|}{n^2}$, $q^*=(\pi,0)$ and we require
$s'>3$, it can be satisfied if  $s>r+3$.

For (\ref{ha}), we have the following lemma.
\begin{Lemma}\label{key}
For $n$ large enough, the Hamiltonian flow generated by $H_n(q,p)$
does not admit any Lagrangian torus with rotation vector
$\omega=(\omega_1,\ldots,\omega_d)$ satisfying
\[|\omega_1|<n^{-\frac{a}{2}-\epsilon},\] where $\epsilon>0$ is independent of $n$.
\end{Lemma}Lemma \ref{key} will be proved with variational method. First
of all, we put it into the Lagrangian formalism. Let
$\sigma_n=n^{-a}$. The Lagrangian function corresponding to
(\ref{Lo}) is
\begin{equation}\label{L1}
\begin{split}
L_n\left(q_1,Q,\dot{q}_1,\dot{Q}\right)=\frac{1}{2}|\dot{Q}|^2+
\frac{1}{2}|\dot{q}_1|^2+\sigma_n(1-\cos(q_1))+v_n(q_1,q_2),
\end{split}\end{equation}
 where $Q=(q_2,\ldots,q_d)$. $L_n(q_1,Q,\dot{q}_1,\dot{Q})$  can be considered as a perturbation coupling of a rotator with
 $d-1$ degrees of freedom and a perturbation with the Lagrangian
 function
 \begin{equation}\label{simpendul_1}
A_n(q_1,\dot{q}_1)=\frac{1}{2}|\dot{q}_1|^2+\sigma_n(1-\cos(q_1)),
\end{equation}
which corresponds to the Hamiltonian via Legendre transformation
\begin{equation}
h_n(q_1,p_1)=\frac{1}{2}|p_1|^2-\sigma_n(1-\cos q_1).
\end{equation}

\subparagraph{The action of the simple pendulum} Each solution of
the Lagrangian equation determined by $A_n$, denoted by $q_1(t)$,
determines an orbit $(q_1(t),p_1(t))$ of the Hamiltonian flow
generated by $h_n$. Each orbit stays in certain energy level set
$(q_1,p_1)\in h_n^{-1}(e)$. Under the boundary condition that
$t_0=0,\ t_1=\pi$ (or $t_1=\pi,t_2=2\pi$), there is a unique
correspondence between $t_1-t_0$ and the energy, denoted by
$e(t_1-t_0)$, such that the determined orbit stays in the energy
level set $h_n^{-1}(e(t_1-t_0))$. More precisely, we have the
following lemma.
\begin{Lemma}\label{energ} Let $\bar{q}_1$ be the solution of $A_n$ on
$(t_0,\bar{t}_1)$ satisfying the boundary conditions
\begin{equation*}
\begin{cases}
\bar{q}_1(t_0)=0,\\
\bar{q}_1(\bar{t}_1)=\pi,
\end{cases}
\end{equation*}  $e(\bar{t}_1-t_0)$ be the energy of $\bar{q}_1$, i.e. $(\bar{q}_1,\bar{p}_1)\in h_n^{-1}(e(\bar{t}_1-t_0))$ and
$\omega_1$ be the average speed of $\bar{q}_1$ on $(t_0,\bar{t}_1)$,
then
\begin{equation}\label{hh_1}
e(\bar{t}_1-t_0)\sim\sigma_n\exp\left(-\frac{C\sqrt{\sigma_n}}{|\omega_1|}\right),
\end{equation}
where $f\sim g$ means that $\frac{1}{C}g<f<Cg$ holds for some
constant $C>0$, $\sigma_n=n^{-a}$. \end{Lemma}

 \Proof By the definition,
we have
\[\frac{1}{2}|\dot{\bar{q}}_1|^2-\sigma_n(1-\cos(\bar{q}_1))=e(\bar{t}_1-t_0),\]
hence
\[|\dot{\bar{q}}_1|=\sqrt{2(e(\bar{t}_1-t_0)+\sigma_n(1-\cos(\bar{q}_1)))}.\]
Since the average speed of $\bar{q}_1$ is $\omega_1$, by a direct
calculation, we have
\[\frac{\pi}{|\omega_1|}=\int_{t_0}^{\bar{t}_1}dt=\int_0^{\pi}\frac{d\bar{q}_1}{\sqrt{2(e(\bar{t}_1-t_0)+\sigma_n(1-\cos(\bar{q}_1)))}}\sim\frac{1}
{\sqrt{\sigma_n}}\ln\left(\frac{\sigma_n}{e(\bar{t}_1-t_0)}\right),\]
moreover,
\[e(\bar{t}_1-t_0)\sim\sigma_n\exp\left(-\frac{C\sqrt{\sigma_n}}{|\omega_1|}\right),\] which complete the proof
of Lemma \ref{energ}.\End

It is easy to see that Lemma \ref{energ} also holds for
\begin{equation*}
\begin{cases}
\bar{q}_1(\bar{t}_1)=\pi,\\
\bar{q}_1(t_2)=2\pi.
\end{cases}
\end{equation*}

The following lemma implies that the actions along orbits in the
neighborhood of the separatix of the  pendulum does not change too
much with respect to a small change in speed (time).

 \begin{Lemma}\label{8_1} Let
$\bar{t}_1, \tilde{t}_1\in [t_0,t_2]$. Let $\bar{q}_1(t)$ be a
solution of $A_n$ on $(t_0,\bar{t}_1)$ and $(\bar{t}_1,t_2)$ with
boundary conditions respectively
 \begin{equation*}\begin{cases}\bar{q}_1(t_0)=0,\\
\bar{q}_1(\bar{t}_1)=\pi,\end{cases}
\begin{cases}\bar{q}_1(\bar{t}_1)=\pi,\\
\bar{q}_1(t_2)=2\pi,\end{cases}\
\end{equation*}
and let $\tilde{q}_1(t)$ be a solution of $A_n$ on
$(t_0,\tilde{t}_1)$ and $(\tilde{t}_1,t_2)$ with boundary conditions
respectively
\begin{equation*}
\begin{cases}\tilde{q}_1(t_0)=0,\\
\tilde{q}_1(\tilde{t}_1)=\pi,\end{cases}
\begin{cases}\tilde{q}_1(\tilde{t}_1)=\pi,\\
\tilde{q}_1(t_2)=2\pi.
\end{cases}
\end{equation*} Let
$\bar{\omega}'_1$ and $\bar{\omega}''_1$ be the average speed of
$\bar{q}_1$ on $(t_0,\bar{t}_1)$ and $(\bar{t}_1,t_2)$ respectively.
Let $\tilde{\omega}'_1$ and $\tilde{\omega}''_1$ be the average
speed of $\tilde{q}_1$ on $(t_0,\tilde{t}_1)$ and
$(\tilde{t}_1,t_2)$ respectively. We set
\[|\omega_1|=\max\left\{|\bar{\omega}'_1|,|\bar{\omega}''_1|,|\tilde{\omega}'_1|,|\tilde{\omega}''_1|\right\},\] then
\begin{equation}\left|\int^{t_2}_{t_0}A_n(\bar{q}_1,\dot{\bar{q}}_1)dt-\int^{t_2}_{t_0}A_n(\tilde{q}_1,\dot{\tilde{q}}_1)dt\right|\leq
C_1|\bar{t}_1-\tilde{t}_1|\sigma_n\exp\left(-\frac{C_2\sqrt{\sigma_n}}{|\omega_1|}\right).\end{equation}
\end{Lemma}

\Proof The proof follows the similar idea of Lemma 4 in \cite{B2}.
Let $q_1(t)$ be a solution of $A_n$ on $(t_0,t_1)$ and $(t_1,t_2)$
with boundary conditions respectively
\begin{equation*}\begin{cases}q_1(t_0)=0,\\
q_1(t_1)=\pi,\end{cases}
\begin{cases}q_1(t_1)=\pi,\\
q_1(t_2)=2\pi.
\end{cases}\end{equation*} We consider the function
\begin{align*}
L(t_1)=&\int_{t_0}^{t_1}A_n(q_1,\dot{q}_1)dt+\int_{t_1}^{t_2}A_n(q_1,\dot{q}_1)dt,\\
=&\int_0^\pi\sqrt{2(e(t_1-t_0)+V(q_1))}dq_1-e(t_1-t_0)(t_1-t_0)\\
&+\int_\pi^{2\pi}\sqrt{2(e(t_2-t_1)+V(q_1))}dq_1-e(t_2-t_1)(t_2-t_1),
\end{align*}
where
\[V(q_1)=\sigma_n(1-\cos(q_1)),\]and $e(\Delta t)$ denotes the energy of the
orbit of the pendulum moving half a turn in time $\Delta t$. The
quantity $e(\Delta t)$ is differentiable with respect to $\Delta t$,
then
\begin{equation}\label{hhh}
\begin{split}
\frac{dL(t_1)}{dt_1}=&\int_0^\pi\frac{\dot{e}(t_1-t_0)}{\sqrt{2(e(t_1-t_0)+V(q_1))}}dq_1-\dot{e}(t_1-t_0)(t_1-t_0)\\
&-e(t_1-t_0)-\int_\pi^{2\pi}\frac{\dot{e}(t_2-t_1)}{\sqrt{2(e(t_2-t_1)+V(q_1))}}dq_1\\
&+\dot{e}(t_2-t_1)(t_2-t_1)+e(t_2-t_1),\\
 =&\int_{t_0}^{t_1}\dot{e}(t_1-t_0)dt-\dot{e}(t_1-t_0)(t_1-t_0)-e(t_1-t_0)\\
 &-\int_{t_1}^{t_2}\dot{e}(t_2-t_1)dt+\dot{e}(t_2-t_1)(t_2-t_1)+e(t_2-t_1),\\
 =&e(t_2-t_1)-e(t_1-t_0).\\
 \end{split}
\end{equation}
Thus, we have \[\left|\frac{dL(t_1)}{dt_1}\right|\leq
|e(t_2-t_1)|+|e(t_1-t_0)|.\]Integrate from $\bar{t}_1$ to
$\tilde{t}_1$ and from (\ref{hh_1}), it follows that
\[\left|\int^{t_2}_{t_0}A_n(\bar{q}_1,\dot{\bar{q}}_1)dt-\int^{t_2}_{t_0}A_n(\tilde{q}_1,\dot{\tilde{q}}_1)dt\right|\leq C_1|\bar{t}_1-\tilde{t}_1|\sigma_n\exp\left(-\frac{C_2\sqrt{\sigma_n}}{|\omega_1|}\right),\]
which completes the proof of Lemma \ref{8_1}. \End

\subparagraph{The velocity of the action minimizing orbit} Once the
function $q_1(t)$ is fixed, the function $Q(t)$ is the solution of
the Euler-Lagrange equation with the non autonomous Lagrangian
\begin{equation}\label{LQ}
\frac{1}{2}|\dot{Q}(t)|^2+v_n\left(q_1(t), q_2(t)\right),
\end{equation}where $Q(t)=(q_2(t),\ldots,q_d(t))$.

\begin{Lemma}\label{9_1} Let $(q_1(t),Q(t))$ be the  orbit of $L_n$ with
rotation vector $\omega$, then for any $t',t''\in \R$ and $t\in
[t',t'']$ we have
\begin{equation}\label{vv}\left|\dot{Q}(t)-\frac{Q(t'')-Q(t')}{t''-t'}\right|\leq
C|\omega_1|^2.\end{equation}\end{Lemma}

\Proof By the Euler-Lagrange equation, we have
\[\ddot{q}_i(t)=0,\quad\text{for}\quad i=3,\ldots,d,\]hence
$\dot{q}_i(t)=\text{const.}$, (\ref{vv}) is verified obviously  for
$q_i(t)$, $i=3,\ldots,d$. We just need to consider $q_2(t)$. Let
$q_1(t_0)=0$ and $q_1(t_2)=2\pi$. It suffices to prove that for
$t\in [t',t'']\subset[t_0,t_2]$
\begin{equation}\label{vv1}\left|\dot{q}_2(t)-\frac{q_2(t'')-q_2(t')}{t''-t'}\right|\leq
C|\omega_1|^2.\end{equation} From the Euler-Lagrange equation,
\[\ddot{q}_2(t)=\frac{\partial v_n}{\partial q_2}(q_1(t),q_2(t)),\]
together with $||v_n||_{C^r}\sim|\omega_1|^{s'}$, we obtain
\[\ddot{q}_2(t)\leq C_1|\omega_1|^{s'}.\] Integrate the two
sides of the inequality above from $t'$ to $t''$, we have
\[|\dot{q}_2(t'')-\dot{q}_2(t')|\leq C_2
|\omega_1|^{s'}|t''-t'|.\] It follows from (\ref{ap}) that
$|t''-t'|\leq C_3|\omega_1|^{-1}$. Hence
\[|\dot{q}_2(t'')-\dot{q}_2(t')|\leq C_4|\omega_1|^{s'-1}.\] Since $s'>3$, we have
\[|\dot{q}_2(t'')-\dot{q}_2(t')|\leq
C_4|\omega_1|^2.\]This completes the proof. \End

\subparagraph{Proof of Lemma \ref{key}} Based on the minimal
property of the orbits on an invariant Lagrangian torus and its
graph property, passing through each $x\in\T^d$, there is a unique
minimal curve $q(t)$ with rotation vector $\omega$ if the
Hamiltonian flow generated by $H_n$ admits a Lagrangian torus with
rotation vector $\omega$. Hence, it is sufficient to prove the
existence of some point in $\T^d$ where no minimal curve passes
through.

Indeed, any minimal curve does not pass through the subspace
$(\pi,0)\times\T^{d-2}$. It implies Lemma \ref{key}. Let
 us assume the contrary, namely, there exists $\bar{t}_1$ such that
\[q_1(\bar{t}_1)=\pi,\quad q_2(\bar{t}_1)=0,\] where
$q(t)=(q_1,q_2,\ldots,q_d)(t)$ is a minimal curve in the universal
covering space $\R^d$. Because of $\omega_1\neq 0$, there exist
$t_0$ and $t_2$ such that
\[q_1(t_0)=0,\quad q_1(t_2)=2\pi.\]
Obviously, $t_0<\bar{t}_1<t_2$ and
\[t_2-t_0\sim\frac{1}{|\omega_1|}.\]

Let $\tilde{t}_1$ be the last time before $\bar{t}_1$ or the first
time after $\bar{t}_1$ such that
\[|q_2(\tilde{t}_1)-q_2(\bar{t}_1)|=\pi.\]
It is easy to see that,
\[|\tilde{t}_1-\bar{t}_1|\sim\frac{1}{|\omega_2|}.\]Since $|\omega_2|\sim 1$, then
\[|\tilde{t}_1-\bar{t}_1|\leq C_0.\]
 Without loss of generality, one can assume
$\omega_1>0$ and $\omega_2>0$. Consider a solution $\tilde{q}_1$ of
$A_n$ on $(t_0,\tilde{t}_1)$ and on $(\tilde{t}_1,t_2)$ with
boundary conditions respectively
\begin{equation*}
\begin{cases}
\tilde{q}_1(t_0)=q_1(t_0)=0,\\
\tilde{q}_1(\tilde{t}_1)=q_1(\bar{t}_1)=\pi,
\end{cases}
\quad
\begin{cases}
\tilde{q}_1(\tilde{t}_1)=q_1(\bar{t}_1)=\pi,\\
\tilde{q}_1(t_2)=q_1(t_2)=2\pi.
\end{cases}
\end{equation*}
Since $q$ is assumed to be a minimal curve, we have
\begin{equation}\label{star}\int_{t_0}^{t_2}L_n(\tilde{q}_1,Q,\dot{\tilde{q}}_1,\dot{Q})dt-\int_{t_0}^{t_2}L_n(q_1,Q,\dot{q}_1,\dot{Q})dt\geq
0.\end{equation} See Fig.1, where
$x_1=(q_1(\bar{t}_1),q_2(\bar{t}_1))=(\pi,0)$,
$x_0=(q_1(t_0),q_2(t_0))=(0,q_2(t_0))$,
$x_2=(q_1(t_2),q_2(t_2))=(2\pi,q_2(t_2))$, $\tilde{x}'_1=(\pi,-\pi)$
and $\tilde{x}''_1=(\pi,\pi)$.

\input{figure1.TpX}

 $(\tilde{q}_1(t),q_2(t))$ passes through the point $\tilde{x}'_1$ or $\tilde{x}''_1$. Thus, by the construction of $L_n$, we
 obtain from (\ref{star}) that
\begin{equation}\label{AandP_1}
\int_{t_0}^{t_2}A_n(\tilde{q}_1,\dot{\tilde{q}}_1)dt-\int_{t_0}^{t_2}A_n(q_1,\dot{q_1})dt\geq
\int_{t_0}^{t_2}v_n(q_1,q_2)dt-\int_{t_0}^{t_2}v_n(\tilde{q}_1,q_2)dt.
\end{equation}
By the definition of $v_n$ as (\ref{v}), we find
\[(\tilde{q}_1(t),q_2(t))\cap\text{supp}\,v_n=\emptyset,\quad\text{for}\quad t\in (t_0,t_2).\]
In fact, if there would exist $\hat{t}$ such that
$(\tilde{q}_1(\hat{t}),q_2(\hat{t}))\in \text{supp}\,v_n$, without
loss of generality, one can assume $\hat{t}>\tilde{t}_1$. By Lemma
\ref{9_1}, for any $t\in [\tilde{t}_1,\hat{t}]$,
\[\dot{q}_2(t)\leq C_1,\]hence,
\[\hat{t}-\tilde{t}_1\geq C_2,\]
where $C_1$, $C_2$ are constants independent of $n$. Consequently
\[|\tilde{q}_1(\hat{t})-\tilde{q}_1(\tilde{t}_1)|\geq
C_3|\omega_1|>R_n,\] where $R_n$ is the radius of the support of
$v_n$. It is impossible.

Hence, we have
\[\int_{t_0}^{t_2}v_n(q_1,q_2)dt-\int_{t_0}^{t_2}v_n(\tilde{q}_1,q_2)dt=\int_{t_0}^{t_2}v_n(q_1,q_2)dt.\]
By the construction of $v_n$ and the minimality of $(q_1,Q)$, a
simple calculation shows
\begin{equation}\label{tt}
\int_{t_0}^{t_2}v_n(q_1,q_2)dt\geq {|\omega_1|}^{\lambda},
\end{equation}where $\lambda$ is a positive constant. Consequently, if follows from (\ref{AandP_1}) that
\begin{equation}\label{contr1_1}
\int_{t_0}^{t_2}A_n(\tilde{q}_1,\dot{\tilde{q}}_1)dt-\int_{t_0}^{t_2}A_n(q_1,\dot{q}_1)dt\geq
 {|\omega_1|}^{\lambda}.
\end{equation}

On the other hand, consider a solution $\bar{q}_1$ of $A_n$ on
$(t_0,\bar{t}_1)$ and on $(\bar{t}_1,t_2)$ with boundary conditions
respectively
\begin{equation*}
\begin{cases}
\bar{q}_1(t_0)=q_1(t_0)=0,\\
\bar{q}_1(\bar{t}_1)=q_1(\bar{t}_1)=\pi,
\end{cases}
\quad
\begin{cases}
\bar{q}_1(\bar{t}_1)=q_1(\bar{t}_1)=\pi,\\
\bar{q}_1(t_2)=q_1(t_2)=2\pi.
\end{cases}
\end{equation*}
Along both of which the action of $A_n$ achieves its minimum. Thus,
we have
\[\int_{t_0}^{t_2}A_n(q_1,\dot{q}_1)dt\geq\int_{t_0}^{t_2}A_n(\bar{q}_1,\dot{\bar{q}}_1)dt.\]

We compare the action
$\int_{t_0}^{t_2}A_n(\tilde{q}_1,\dot{\tilde{q}}_1)dt$ with the
action $\int_{t_0}^{t_2}A_n(\bar{q}_1,\dot{\bar{q}}_1)dt$  in the
alternative cases, which is  based on the choices of $\tilde{t}_1$.
See Fig.2, where $\bar{t}=\frac{t_0+t_2}{2}$.

\input{figure2.TpX}

\noindent Case 1: $|\bar{t}_1-\bar{t}|\leq C_0$.

In this case, the average speed of $\bar{q}_1$ on $(t_0,\bar{t}_1)$
and $(\bar{t}_1,t_2)$ have the same quantity order as $|\omega_1|$.
By $|\tilde{t}_1-\bar{t}_1|\leq C_0$, we have
$|\tilde{t}_1-\bar{t}|\leq 2C_0$. Hence the average speed of
$\tilde{q}_1$ on $(t_0,\tilde{t}_1)$ and $(\tilde{t}_1,t_2)$ have
also the same quantity order as $|\omega_1|$. Thus, Lemma \ref{8_1}
implies
\[\int_{t_0}^{t_2}A_n(\tilde{q}_1,\dot{\tilde{q}}_1)dt-\int_{t_0}^{t_2}A_n(\bar{q}_1,\dot{\bar{q}}_1)dt\leq
C_4\sigma_n\exp\left(-\frac{C_5\sqrt{\sigma_n}}{|\omega_1|}\right).\]

\noindent Case 2: $|\bar{t}_1-\bar{t}|> C_0$.

In this case, we take $\tilde{t}_1$ such that
$|\tilde{t}_1-\bar{t}|\leq |\bar{t}_1-\bar{t}|$, which can be
achieved by the suitable choice of the position of
$\tilde{q}_1(\tilde{t}_1)$. More precisely,
\begin{itemize}
  \item  if $\bar{t}_1>\bar{t}+C_0$ (Case 2a in Fig.2), we choose $\tilde{t}_1$ as the
last time before $\bar{t}_1$, i.e. $(\tilde{q}_1(\tilde{t}_1),
q_2(\tilde{t}_1))=\tilde{x}'_1$ in Fig.1;
  \item  if
$\bar{t}_1<\bar{t}-C_0$ (Case 2b in Fig.2), we choose $\tilde{t}_1$
as the first time after $\bar{t}_1$, i.e.
$(\tilde{q}_1(\tilde{t}_1), q_2(\tilde{t}_1))=\tilde{x}''_1$ in
Fig.1.
\end{itemize}
  For both cases 2a and 2b, it
follows from (\ref{hhh}) that
\[\int_{t_0}^{t_2}A_n(\tilde{q}_1,\dot{\tilde{q}}_1)dt-\int_{t_0}^{t_2}A_n(\bar{q}_1,\dot{\bar{q}}_1)dt\leq
0.\] Hence, for any $\bar{t}_1\in (t_0,t_2)$, we can find
$\tilde{t}_1$ such that
\begin{align*}
\int_{t_0}^{t_2}A_n(\tilde{q}_1,\dot{\tilde{q}}_1)dt-\int_{t_0}^{t_2}A_n(q_1,\dot{q}_1)dt&\leq
\int_{t_0}^{t_2}A_n(\tilde{q}_1,\dot{\tilde{q}}_1)dt-\int_{t_0}^{t_2}A_n(\bar{q}_1,\dot{\bar{q}}_1)dt,\\
&\leq
C_4\sigma_n\exp\left(-\frac{C_5\sqrt{\sigma_n}}{|\omega_1|}\right).
\end{align*}
Since \[|\omega_1|\leq n^{-\frac{a}{2}-\epsilon}.\] It is easy to
see that for $n$ large enough,
\[C_4\sigma_n\exp\left(-\frac{C_5\sqrt{\sigma_n}}{|\omega_1|}\right)\leq
{|\omega_1|}^{\lambda},\]where $\sigma_n=n^{-a}$, which contradicts
to (\ref{contr1_1}) for large $n$. This completes the proof of Lemma
\ref{key}.\End

\paragraph{Destruction of Lagrangian torus with an arbitrary rotation
vector} By (\ref{ap}), for every non resonant rotation vector
$\omega=(\omega_1,\ldots,\omega_d)$ $(d\geq 2)$, there exists a
sequence of integer vector $k_n\in\Z^d$ satisfying $|k_n|\rightarrow
\infty$ as $n\rightarrow \infty$ such that
\[|\langle k_n, \omega\rangle|<\frac{C}{|k_n|^{d-1}}.\]
\subparagraph{Transformation of coordinates} We choose a sequence of
$k_n\in\Z^d$ satisfying (\ref{ap}) and an integer vector sequence
$k'_n$ such that $\langle k'_n,k_n\rangle=0$. In addition, select
$d-2$ integer vectors $l_{n3},\ldots,l_{nd}$ such that
$k_n,k'_n,l_{n3},\ldots,l_{nd}$ are pairwise orthogonal. Let
\begin{equation}\label{K}
K_n=(k_n,k'_n,l_{n3},\ldots,l_{nd})^t. \end{equation} We choose the
transformation of the coordinates
\[q=K_nx.\]
Let $p$ denotes the dual coordinate of $q$ in the sense of  Legendre
transformation, i.e. $p=\frac{\partial L}{\partial\dot{q}}$, it
follows that
\[y=K^t_np,\]where $K^t_n$ denotes the transpose of $K_n$. We
set
\[\Phi_n=\begin{pmatrix}
K_n&\ \\
\ &K^{-t}_n \end{pmatrix},\]then
\[\begin{pmatrix}q\\p\end{pmatrix}=\Phi_n\begin{pmatrix}x\\y\end{pmatrix}.\]
It is easy to verify that \[\Phi_n^tJ_0\Phi_n=J_0,\] where
\[J_0=\begin{pmatrix}
\textbf{0}&\textbf{1}\\
-\textbf{1}&\textbf{0}\end{pmatrix},\] where $\textbf{1}$ denotes a
$d\times d$ unit matrix. Hence, $\Phi_n$ is a symplectic
transformation in the phase space.

\begin{Lemma}\label{key1}
If the Hamiltonian flow generated by $\tilde{H}_n(x,y)$ admits a
Lagrangian torus with rotation vector $\omega$, then the Hamiltonian
flow generated by $H_n(q,p)$ also admits a Lagrangian torus with
rotation vector $K_n\omega$, where $(q,p)^t=\Phi_n(x,y)^t$.
\end{Lemma}

\Proof Let $\tilde{\mathcal {T}}^d$ be the Lagrangian torus admitted
by $\tilde{H}_n(x,y)$, a symplectic form $\Omega$ vanishes on
$\text{T}_x\tilde{\mathcal {T}}^d$ for every $x\in \tilde{\mathcal
{T}}^d$. Since $K_n$ consists of integer vectors, then $\mathcal
{T}^d:=K_n\tilde{\mathcal {T}}^d$ is still a torus. $\Phi_n$ is a
symplectic transformation, hence $\mathcal {T}^d$ is a Lagrangian
torus. From Definition \ref{dd1}, each orbit on $\tilde{\mathcal
{T}}^d$ has the same rotation vector $\omega$. Let
$\tilde{\gamma}(t)$ be a lift of an orbit on $\tilde{\mathcal
{T}}^d$, it follows that
\[\omega=\lim_{t\rightarrow\infty}\frac{\tilde{\gamma}(t)-\tilde{\gamma}(-t)}{2t}.\]
For $\gamma(t)=K_n\tilde{\gamma}(t)$, we have
\begin{align*}
\lim_{t\rightarrow\infty}\frac{\gamma(t)-\gamma(-t)}{2t}&=\lim_{t\rightarrow\infty}\frac{K_n\tilde{\gamma}(t)-K_n\tilde{\gamma}(-t)}{2t};\\
&=K_n\lim_{t\rightarrow\infty}\frac{\tilde{\gamma}(t)-\tilde{\gamma}(-t)}{2t};\\
&=K_n\omega.
\end{align*}
This completes the proof.\End
\begin{Remark}
For $\tilde{\mathcal {T}}^d$ and $\mathcal {T}^d$ in the proof of
Lemma \ref{key1}, we have
\[\text{Vol}\,(\mathcal {T}^d)=|\det K_n|\text{Vol}\,(\tilde{\mathcal
{T}}^d),\]where Vol($\cdot$) denotes the volume of ($\cdot$).
\end{Remark}
\subparagraph{Proof of Theorem \ref{maint}} We construct
$\tilde{H}_n(x,y)$ as follow:
\begin{equation}\label{H}
\tilde{H}_n(x,y)=\frac{1}{2}|y|^2-\tilde{P}_n(x),
\end{equation}
where
\[\tilde{P}_n(x)=\frac{1}{|k_n|^{a+2}}(1-\cos\langle
k_n,x\rangle)+\frac{1}{|k_n|^2}v_n\left(\langle k_n,x\rangle,\langle
k'_n,x\rangle\right),\] where $k'_n$ is the second row of $K_n$ and
$v_n$ is defined by (\ref{v}).  Let $q=K_nx$. In particular, we have
\be\label{cortranfor}
\begin{cases}q_1=\langle k_n,x\rangle,\\
q_2=\langle k'_n,x\rangle.
\end{cases}\ee
By the transformation of coordinates and the Legendre
transformation, the Lagrangian function corresponding to (\ref{H})
is
\begin{equation}\label{ll}
\begin{split}
L_n\left(q_1,Q,\dot{q}_1,\dot{Q}\right)=&\frac{1}{2}\sum^d_{i=3}\frac{|\dot{q}_i|^2}{|l_{ni}|^2}+\frac{|\dot{q}_2|^2}{2|k'_n|^2}\\
&+\frac{1}{|k_n|^2}\left(\frac{1}{2}|\dot{q}_1|^2+\frac{1}{|k_n|^a}(1-\cos(q_1))+v_n(q_1,q_2)\right),
\end{split}
\end{equation} where $Q=(q_2,\ldots,q_d)$.

For the Hamiltonian flow generated by (\ref{H}), by Lemma
\ref{key1}, for  the destruction of the Lagrangian torus
$\tilde{\mathcal {T}}^d$ with rotation vector $\omega$, it is
sufficient to prove that the Euler-Lagrange flow generated by
(\ref{ll}) admits no the Lagrangian torus $\mathcal
{T}^d:=K_n\tilde{\mathcal {T}}^d$ with rotation vector $K_n\omega$.
 Let $K_n\omega=(\omega_1,\omega_2,\ldots,\omega_d)$.

It is easy to see that there exists an integer vector $k'_n$ such
that
\begin{equation}\label{kn1}
\langle k_n,k'_n\rangle=0\quad\text{and}\quad |\langle
k'_n,\omega\rangle|\sim 1,
\end{equation}
i.e. $\omega_2\sim 1$. In fact, it suffices to consider
$k'_n\in\Z^3$. Let $k'_n=(k'_{n1},k'_{n2},k'_{n3})$, then for
$k'_n\in\Z^d$, one can take
$k'_n=(k'_{n1},k'_{n2},k'_{n3},0,\ldots,0)$ to verify (\ref{kn1}).

Since $\omega$ is non-resonant, then  $|k'_n|\rightarrow\infty$, for
$n\rightarrow \infty$. Let $\theta$ be the angle between $k'_n$ and
$\omega$, then
 \[|\langle k'_n,\omega\rangle|=|k'_n||\omega||\cos\theta|,\]
 where $|k'_n|$ is determined by (\ref{de}) below.
 Let $\Pi$ be the plane
orthogonal with respect to $k_n$. Let $S_{R\alpha}\subset\Pi$ be the
sector with cental point $(0,0,0)$, central angle $\alpha$ and
radius $R$ satisfying $\alpha=\frac{C_1}{|k'_n|}$ and the angle
between $\omega$ and one of the radii is equal to
$\frac{\pi}{2}-\frac{C_2}{|k'_n|}$, where $C_2>C_1$. To achieve
(\ref{kn1}), it is sufficient to find an integer point
$m=(m_1,m_2,m_3)\in\Z^3$ satisfying
\begin{equation}\label{mm}
|m|\sim |k'_n|\quad\text{and}\quad m\in S_{R\alpha}.
\end{equation}
In deed, we take $k'_n=(m_1-0,m_2-0,m_3-0)$. Since
\[\cos\theta\sim\cos(\frac{\pi}{2}-\frac{1}{|k'_n|})=\sin\frac{1}{|k'_n|}\sim \frac{1}{|k'_n|}, \]
we have $|\langle k'_n,\omega\rangle|\sim 1$. It is easy to see that
there exists a suitable constant $r(|k'_n|)$ only depending on
$|k'_n|$ such that the square of area $(r(|k'_n|))^2$ contains at
least one integer point. We take
\begin{equation}\label{de}
|k'_n|\sim (r(|k'_n|))^\kappa, \end{equation}where $\kappa\gg 1$,
then it can be concluded that the integer satisfying (\ref{mm}) does
exist.

Replacing $n$ by $|k_n|$ in the proof of Lemma \ref{key}, we have
that the Euler-Lagrange flow generated by (\ref{ll}) does not admit
any Lagrangian torus with rotation vector satisfying
\[|\omega_1|<|k_n|^{-\frac{a}{2}-\epsilon}.\]
From the construction of $K_n$, $|\omega_1|=|\langle k_n,
\omega\rangle|$ which together with (\ref{ap}) implies
\[|\omega_1|<\frac{C}{|k_n|^{d-1}}.\]
Based on Lemma \ref{key1},  it suffices to take
\[\frac{C}{|k_n|^{d-1}}\leq |k_n|^{-\frac{a}{2}-\epsilon},\] which
implies
\begin{equation}\label{a}
a<2d-2-2\epsilon.
\end{equation}

It follows from (\ref{v}) and (\ref{H}) that
\begin{align*}
||\tilde{H}_n&(x,y)-H_0(y)||_{C^r}\\
&=||\tilde{P}_n(x)||_{C^r},\\
&=|k_n|^{-2}\left(|k_n|^{-a}||1-\cos\langle k_n,
x\rangle||_{C^r}+||v_n(\langle k_n,x\rangle, \langle k'_n,
x\rangle)||_{C^r}\right),\\
&\leq |k_n|^{-2}\left(C_1|k_n|^{-a+r}+C_2|k_n|^{-s'(d-1)+r}\right),\\
&\leq C_3\left(|k_n|^{r-a-2}+|k_n|^{r-3(d-1)-2}\right).
\end{align*}

To complete the proof of Theorem \ref{maint}, it is enough to make
$r-a-2<0$ and $r-3(d-1)-2<0$, which together with (\ref{a}) implies
\[r<2d-2\epsilon.\]Taking $\delta=3\epsilon$, this completes the proof of
Theorem \ref{maint}.\End

\subsubsection{$C^\omega$ case}
We will prove the following theorem:
\begin{Theorem}\label{maint11} Given an integrable Hamiltonian $H_0$ with $d\ (d\geq 3)$ degrees of freedom,
 a rotation vector $\omega$ and a small positive constant
$\delta$, there exists a sequence of $C^\omega$ Hamiltonians
$\{H_n\}_{n\in \N}$ such that $H_n\rightarrow H_0$ in
$C^{d+1-\delta}$ topology  and the Hamiltonian flow generated by
$H_n$ admits no  Lagrangian torus with the rotation vector $\omega$.
\end{Theorem}
\paragraph{Construction of $H_n$}
$P_n(x)$ is constructed as follow:
\begin{equation*}
P_n(x)=\frac{1}{|k_n|^{d+1-\epsilon}}(1-\cos\langle
k_n,x\rangle)+\mu_n \frac{1}{|k_n|^{d+1-\epsilon}}(1-\cos\langle
k_n,x\rangle)\cos\langle k'_n,x\rangle,
\end{equation*} where $k_n$, $k'_n$ are the first two rows of $K_n$ defined as (\ref{K}), $\epsilon$ is a given small positive constant and $\mu_n$ satisfies
\begin{equation}\label{mu}
\mu_n\sim
\exp\left(-|k_n|^{\frac{d}{2}-\frac{1}{2}+\frac{\epsilon}{3}}\right).
\end{equation}
A simple calculation implies that for $\delta=3\epsilon$
\[||H_n(x,y)-H_0(y)||_{C^{d+1-\delta}}=||P_n(x)||_{C^{d+1-\delta}}\rightarrow 0\quad \text{as}\ n\rightarrow \infty.\]
From the transformation matrix of coordinates (\ref{K}), let
$\delta_n=\frac{1}{|k_n|^{d-1-\epsilon}}$, the Lagrangian function
corresponding to (\ref{Lo}) is
\begin{equation}\label{L}
\begin{split}
L_n\left(q_1,Q,\dot{q}_1,\dot{Q}\right)&=\frac{1}{2}\sum^d_{i=3}\frac{|\dot{q}_i|^2}{|l_{ni}|^2}+\frac{|\dot{q}_2|^2}{2|k'_n|^2}+\frac{1}{|k_n|^2}\left(\frac{1}{2}|\dot{q}_1|^2+\delta_n(1-\cos(q_1))\right)\\
&\quad+\frac{1}{|k_n|^2}\left(\mu_n\delta_n(1-\cos(q_1))\cos(q_2)\right),
\end{split}\end{equation} where $Q=(q_2,\ldots,q_d)$.

$L_n(q_1,Q,\dot{q}_1,\dot{Q})$  can be considered as a perturbation
coupling of a rotator with
 $d-1$ degrees of freedom and a perturbation with the Lagrangian
 function \begin{equation}\label{simpendul}
A_n(q_1,\dot{q}_1)=\frac{1}{2}|\dot{q}_1|^2+\delta_n(1-\cos(q_1)),
\end{equation} which corresponds to the Hamiltonian via Legendre transformation
\begin{equation}
h_n(q_1,p_1)=\frac{1}{2}|p_1|^2-\delta_n(1-\cos q_1).
\end{equation}
We denote the coupling perturbation by
\begin{equation}\label{per}
\tilde{P}_n(q_1,Q)=\mu_n\delta_n(1-\cos(q_1))\cos(q_2).
\end{equation}

\paragraph{Melnikov function} In the following, we give some approximation
lemmas on the actions of (\ref{simpendul}) and (\ref{per}) along the
minimal orbits of $L_n$ by the calculation of Melnikov function.

For $t\in \R$, the separatrix $\hat{q}_1(t)$ of $A_n$ as
(\ref{simpendul}) satisfying $\dot{\hat{q}}_1>0$ and
$\hat{q}_1(0)=\pi$ is \begin{equation}\label{separax}\begin{cases}
\hat{q}_1(t)=4\arctan\left(\exp(\sqrt{\delta_n}t)\right),\\
\dot{\hat{q}}_1(t)=\frac{2\sqrt{\delta_n}}{\cosh(\sqrt{\delta_n}t)}.
\end{cases}
\end{equation} Let the separatrix of $A_n$ takes value $\pi$ at
$t_1$ and $Q(t_1)=Q_1$, then for a given rotation vector $\omega$,
we define the Melnikov function as
\begin{equation*} M_n(\omega, Q_1,t_1)=\delta_n\int_\R
\left(1-\cos(\hat{q}_1(t-t_1))\right)\cos\left(\langle
k'_n,\omega(t-t_1)\rangle+q_2(t_1)\right) dt. \end{equation*}
Namely, $M_n$ is the integral of the coupling perturbation
(\ref{per}) along the separatrix of $A_n$. It can be explicitly
calculated as follow

\begin{equation}\label{mel} M_n(\omega, Q_1,t_1)=2\pi\frac{\langle \omega,
k'_n\rangle}{\sinh(\frac{\pi\langle \omega,
k'_n\rangle}{2\sqrt{\delta_n}})}\cos(q_2(t_1)). \end{equation} It is
easy to see that $M_n$ only depends on $\omega_2=\langle
\omega,k'_n\rangle$ and $q_2(t_1)$ on which $M_n$ is $2\pi$
periodic. For the simplicity of notations, we denote $M_n(\omega,
Q_1, t_1)$ by
\begin{equation}\label{Mn}
 M_n(\omega_2, Q_1, t_1)=2\pi\frac{\omega_2}{\sinh(\frac{\pi\omega_2}{2\sqrt{\delta_n}})}\cos(q_2(t_1)).
\end{equation}
Next, we work on the universal covering space of $\T^d$. By
(\ref{Mn}), a simple calculation implies the next lemma
\begin{Lemma}\label{10} If $\bar{q}_2(t')\,\text{mod}\,2\pi=0$ and
$\tilde{q}_2(t'')\,\text{mod}\,2\pi=\pi$. Let
\[Q'=(\bar{q}_2(t'),q_3(t'),\ldots,q_d(t'))\quad\text{and}\quad
Q''=\tilde{q}_2(t''),q_3(t''),\ldots,q_d(t'')),\] then for $n$ large
enough
\begin{equation}\label{mm} M_n(\omega_2, Q',t')-M_n(\omega_2, Q'',t'')\geq
\exp\left(-\frac{\lambda}{\sqrt{\delta_n}}\right),\end{equation}
where $\lambda$ is a positive constant independent of  $n$.
\end{Lemma}

\Proof By (\ref{kn1}), we have
\[|\omega_2|\sim 1.\] A simple calculation gives
\[ M_n(\omega_2,
Q',t')-M_n(\omega_2, Q'',t'')\geq
\exp\left(-\frac{C_2|\omega_2|}{\sqrt{\delta_n}}\right)\geq\exp\left(-\frac{\lambda}{\sqrt{\delta_n}}\right),\]
where $\lambda$ is a positive constant independent of $n$. This
completes the proof of (\ref{mm}).\End

\paragraph{An approximation lemma} We use $M_n$ to approximate the action
of the perturbation (\ref{per}) along the minimal orbits of $L_n$.

\blm\label{key}Let $(q_1(t),Q(t))$ be the minimal orbit of $L_n$
satisfying $q_1(t_1)=\pi$ with rotation vector $\omega$, then
\begin{description}
\item [(i)]there exist $\tau>0$ and $t_0,\ t_2$ satisfying
\begin{equation}
t_0\leq t_1-\frac{\tau}{\delta'_n}<t_1+\frac{\tau}{\delta'_n}\leq
t_2,
\end{equation}
such that
\[q_1(t_0)=0,\quad q_1(t_2)=2\pi.\] where
\[\delta'_n=\frac{1}{|k_n|^{d-1-\frac{\epsilon}{4}}};\]
\item [(ii)]let
 \begin{equation}\label{mean}
\bar{\omega}=\left(\omega_1,\frac{Q(t_2)-Q(t_0)}{t_2-t_0}\right),
\end{equation}  then
\begin{align*}
\bigg|\delta_n\int^{t_2}_{t_0}\left(1-\cos(q_1(t))\right)\cos(q_2(t))dt&-M_n(\bar{\omega}_2,Q(t_1),t_1)\bigg|\\
&\leq
C\sqrt{\delta_n}\exp\left(-\frac{\lambda}{\sqrt{\delta_n}}\right),
\end{align*}where $C$ is a positive constant independent of $n$ and $\lambda$ is the same as the one in (\ref{mm}).
\end{description}\elm

\Proof The proof of Lemma \ref{key} follows from the ideas of
\cite{B2}. We will prove (i) and (ii) respectively in the following.
For the simplicity of notations, we will use $u\preceq v$ (resp.
$u\succeq v$) to denote $u\leq Cv$ (resp. $u\geq Cv$) for some
positive constant $C$.

\subparagraph{Proof of (i)}We set $\tau=\frac{\pi}{C_0}$, where
$C_0$ is the constant in (\ref{ap}). From (\ref{ap}), it follows
that
\[\frac{2\tau}{\delta'_n}<\frac{2\pi}{|\omega_1|}.\] Hence, we have
either \[t_0\leq t_1-\frac{\tau}{\delta'_n}\quad\text{or}\quad
t_1+\frac{\tau}{\delta'_n}\leq t_2.\] Without loss of generality, we
assume $t_1+\frac{\tau}{\delta'_n}\leq t_2$ and prove $t_0\leq
t_1-\frac{\tau}{\delta'_n}$ in the following.

Let $t_{-1}$ and $t_3$ be the last time to the left of $t_1$ such
that $q_1(t_{-1})=-\pi$ and the first time to the right of $t_1$
such that $q_1(t_3)=3\pi$ respectively. Consider the solution
$\bar{q}$ of $A_n$ on $(t_{-1},t_1)$ satisfying the boundary
condition
\begin{equation*}
\begin{cases}
\bar{q}(t_{-1})=q_1(t_{-1})=-\pi,\\
\bar{q}(t_1)=q_1(t_1)=\pi.
\end{cases}
\end{equation*}
We denote the energy of $\bar{q}$ by $\bar{e}$. Since
$t_1-t_{-1}\succeq |k_n|^{d-1}$, by the deduction as similar as the
one in Lemma \ref{energ}, we have
\begin{equation}\label{hh}
0<\bar{e}\leq
\delta_n\exp\left(-\frac{C|k_n|^{\epsilon}}{\sqrt{\delta_n}}\right).
\end{equation}
We set
\[e(t)=\frac{1}{2}|\dot{q}_1(t)|^2-\delta_n(1-\cos(q_1(t))),\] it is easy to
see that there exists $\bar{t}\in [t_{-1},t_1]$ such that
$e(\bar{t})=\bar{e}$. Indeed, without loss of generality, we assume
by contradiction that $q_1(t)$ lie above $\bar{q}(t)$ in the phase
plane for all $t\in [t_{-1},t_1]$, namely
$\dot{q}_1(t)>\dot{\bar{q}}(t)$ for all $t\in [t_{-1},t_1]$, which
is contradicted by the boundary conditions
$\bar{q}(t_{-1})=q_1(t_{-1})$ and $\bar{q}(t_1)=q_1(t_1)$. Hence,
there exists $\bar{t}\in [t_{-1},t_1]$ such that
$q_1(\bar{t})=\bar{q}(\bar{t})$ and
$\dot{q}_1(\bar{t})=\dot{\bar{q}}(\bar{t})$, moreover we have
$e(\bar{t})=\bar{e}$.

By Euler-Lagrange equation, we can estimate $\dot{e}(t)$. More
precisely,
\[\dot{e}(t)\preceq\mu_n\delta_n\dot{q}_1(t).\]Hence,
integrating from $\bar{t}$ to $t$, we have
\begin{equation}\label{h}
\sup_{t\in [t_{-1},t_1]}|e(t)-\bar{e}|\leq \gamma\mu_n\delta_n,
\end{equation}
where $\gamma$ is a positive constant independent of $n$.

We proceed the proof of (i) by the following three steps.\\
\noindent \textbf{a)} $\dot{q}_1(t_1)>0$.

we assume by contradiction that $\dot{q}_1(t_1)\leq 0$.
$q_1(t_1)=\pi$
 together with (\ref{h}) implies that for $n$ large enough,
 $\dot{q}_1(t_1)<0$. Let us denote by $(t_1-\Delta t,t_1)$ the
 maximal interval on the left of $t_1$ on which $q_1(t)\geq \pi$.
 Since $\dot{q}_1(t_1)<0$, then $\Delta t>0$. From
 $q_1(t_{-1})=-\pi<\pi=q_1(t_1)$, it follows that $t_1-\Delta t>t_{-1}$. Let
 \begin{equation*}
 \tilde{q}(t)=\begin{cases}
 q_1(t)\qquad\quad\  t\in [t_{-1},t_1-\Delta t),\\
 2\pi-q_1(t)\quad t\in [t_1-\Delta t,t_1].
 \end{cases}
\end{equation*}
It is easy to see that $(\tilde{q},Q)$ is still an action minimizing
orbit on $[t_{-1},t_1]$. By the definition of $\Delta t$, we have
that $\tilde{q}(t_1-\Delta t)=\pi$ and
\[\dot{\tilde{q}}((t_1-\Delta t)-)\cdot\dot{\tilde{q}}((t_1-\Delta t)+)\leq 0,\]
where $(t_1-\Delta t)-$ denotes $t$ tends to $t_1-\Delta t$ from the
left side and $(t_1-\Delta t)+$ denotes $t$ tends to $t_1-\Delta t$
from the right side. By Euler-Lagrange equation, we have that
$\dot{\tilde{q}}(t)$ is continuous for $t\in [t_{-1},t_1]$. Hence,
$\dot{\tilde{q}}(t_1-\Delta t)=0$. On the other hand, from
(\ref{hh}) and (\ref{h}), it follows that for $n$ large enough,
$\dot{\tilde{q}}(t_1-\Delta t)\neq 0$. Therefore, we have
$\dot{q}_1(t_1)>0$.\vspace{1em}

\noindent \textbf{b)} $\dot{q}_1(t)>0$ for $t\in \left[
t_1-\frac{\tau}{\delta'_n},t_1\right]$.

 Let $(\tilde{t},t_1]$ be the
maximal interval to the left of $t_1$ on which $\dot{q}_1>0$, hence
we can denote the inverse function by $t(q_1)$. Let
$\Gamma=\gamma\mu_n\delta_n$. It follows from (\ref{h}) that for
$q_1<\pi$,
\begin{align*}t_1+\int_\pi^{q_1}\frac{ds}{\sqrt{2(\bar{e}-\Gamma+\delta_n(1-\cos(s)))}}&\leq
t(q_1),\\
&\leq
t_1+\int_\pi^{q_1}\frac{ds}{\sqrt{2(\bar{e}+\Gamma+\delta_n(1-\cos(s)))}}.\end{align*}

On the other hand, denoting the inverse function of
$\hat{q}_1(t-t_1)$ by $\hat{t}(q_1)$, we have
\[\hat{t}(q_1)=t_1+\int_\pi^{q_1}\frac{ds}{\sqrt{2(\delta_n(1-\cos(s)))}}.\]
Moreover, a simple calculation implies
\begin{equation*}
|t(q_1)-\hat{t}(q_1)|\preceq\frac{\bar{e}+\Gamma}{(\delta_n(1-\cos(q_1)))^{\frac{3}{2}}},
\end{equation*}
which together with (\ref{mu}) implies that
\begin{equation}\label{tt}
|t(q_1)-\hat{t}(q_1)|\preceq\frac{\Gamma}{(\delta_n(1-\cos(q_1)))^{\frac{3}{2}}}.
\end{equation}
 It is easy to see that
\begin{equation}\label{qv}
|\dot{q}_1|+|\dot{\hat{q}}_1|\preceq\sqrt{\delta_n}.
\end{equation}
Let
\[F(t)=\min\{1-\cos(q_1(t)),1-\cos(\hat{q}_1(t-t_1))\}.\]
It follows from (\ref{tt}) and (\ref{qv}) that for $t\in
[\tilde{t},t_2]$
\begin{equation}\label{qq}|q_1(t)-\hat{q}_1(t-t_1)|\preceq
\frac{\Gamma\sqrt{\delta_n}}{(\delta_nF(t))^{\frac{3}{2}}}.\end{equation}
Hence, we have \[\tilde{t}\leq t_1-\frac{\tau}{\delta'_n}.\] In
fact, we assume by contradiction that $\tilde{t}>
t_1-\frac{\tau}{\delta'_n}$. It follows from (\ref{qq}) that for $n$
large enough,
\begin{equation}\label{qhat}
q_1(\tilde{t})\geq
\frac{1}{2}\hat{q}_1\left(-\frac{\tau}{\delta'_n}\right).
\end{equation} If there exists $t^*\in [\tilde{t},t_1]$ such that
$\dot{q}_1(t^*)=0$, then
\[|e(t^*)|\geq\delta_n\exp\left(-\frac{C|k_n|^{\frac{3\epsilon}{4}}}{\sqrt{\delta_n}}\right),\]
which is contradicted by (\ref{h}) and (\ref{mu}). Since
$\dot{q}_1(\tilde{t})\geq 0$, we get $\dot{q}_1(\tilde{t})> 0$,
contradicted by the maximality of $[\tilde{t},t_1]$. Therefore, we
have $\tilde{t}\leq t_1-\frac{\tau}{\delta'_n}$. i.e.
$\dot{q}_1(t)>0$ for $t\in \left[
t_1-\frac{\tau}{\delta'_n},t_1\right]$. \vspace{1em}

\noindent \textbf{c)} $q_1(t_1-\frac{\tau}{\delta'_n})>0$.

From (\ref{qq}) and (\ref{qhat}), we have that for
$t\in\left[t_1-\frac{\tau}{\delta'_n},t_1\right]$
\begin{equation}\label{qq1}
|q_1(t)-\hat{q}_1(t-t_1)|\preceq\Gamma\delta^{-1}_n\exp\left(\frac{C|k_n|^{\frac{3\epsilon}{4}}}{\sqrt{\delta_n}}\right).
\end{equation}
In terms of the definition of $\hat{q}_1(t)$, it follows from
(\ref{mu}) that for $n$ large enough,
\begin{equation}\label{qqq}\frac{1}{2}\hat{q}_1\left(-\frac{\tau}{\delta'_n}\right)\leq
q_1\left(t_1-\frac{\tau}{\delta'_n}\right)\leq
\frac{3}{2}\hat{q}_1\left(-\frac{\tau}{\delta'_n}\right).\end{equation}
Hence, $q_1(t_1-\frac{\tau}{\delta'_n})>0$, which together with
$\dot{q}_1(t)>0$ on $\left[t_1-\frac{\tau}{\delta'_n},t_1\right]$
implies $t_0\leq t_1-\frac{\tau}{\delta'_n}$.

Similarly, we have $t_1+\frac{\tau}{\delta'_n}\leq t_2$. The proof
of $(i)$ is completed.

\subparagraph{Proof of (ii)}
 We let
$\Omega=[t_1-\frac{\tau}{\delta'_n},t_1+\frac{\tau}{\delta'_n}]$.
Since
\begin{align*}
\bigg|\delta_n\int^{t_2}_{t_0}(1-&\cos(q_1(t)))\cos(q_2(t)) dt-M_n(\bar{\omega}_2,Q(t_1),t_1)\bigg|\\
\leq&
C_1\delta_n\left(\int_\Omega|\cos(q_2(t))-\cos\left(\hat{q}_2(t)+q_2(t_1)-\hat{q}_2(t_1)\right)|dt\right.\\
&\left.+\int_\Omega|\cos(q_1(t))-\cos(\hat{q}_1(t-t_1))|dt+\int_{[t_0,t_2]\backslash\Omega}|1-\cos(q_1(t))|dt\right.\\
&\left.+\int_{\R\backslash\Omega}|1-\cos(\hat{q}_1(t-t_1))|dt\right),
\end{align*}
where $\hat{q}_2(t)=q_2(t_0)+\omega_2(t-t_0)$. Hence, it suffices to
prove the following estimates
\begin{equation}\label{1}
\delta_n\int_\Omega|\cos(q_2(t))-\cos\left(\hat{q}_2(t)+q_2(t_1)-\hat{q}_2(t_1)\right)|dt\preceq\epsilon_n,
\end{equation}
\begin{equation}\label{2}
\delta_n\int_\Omega|\cos(q_1(t))-\cos(\hat{q}_1(t-t_1))|dt\preceq\epsilon_n,
\end{equation}
\begin{equation}\label{3}
\delta_n\int_{[t_0,t_2]\backslash
\Omega}|1-\cos(q_1(t))|dt\preceq\epsilon_n,
\end{equation}
\begin{equation}\label{4}
\delta_n\int_{\R\backslash
\Omega}|1-\cos(\hat{q}_1(t-t_1))|dt\preceq\epsilon_n,
\end{equation}
where we set
$\epsilon_n=\sqrt{\delta_n}\exp(-\frac{\lambda}{\sqrt{\delta_n}})$.
We will prove (\ref{1})-(\ref{4}) in the following.

First of all, we prove (\ref{1}).
We set
\begin{equation*}
\bar{q}(t)=\left\{\begin{array}{ll}
\hspace{-0.4em}\hat{q}_1(t-t_1)\phi\left(t-(t_1-\frac{\tau}{\delta'_n})\right),&t_0\leq
t\leq t_1,\\
\hspace{-0.4em}\left(\hat{q}_1(t-t_1)-2\pi\right)\phi\left(-t+(t_1+\frac{\tau}{\delta'_n})\right)+2\pi,&
t_1\leq t\leq t_2,
\end{array}\right.
\end{equation*}
where $\phi$ is a $C^\infty$ function as follow
\begin{equation*}
\phi(t)=\left\{\begin{array}{ll} \hspace{-0.4em}0,& t\leq 0,\\
\hspace{-0.4em}1,& t\geq 1.
\end{array}\right.
\end{equation*}
We set $\hat{Q}(t)=Q(t_0)+\bar{\omega}(t-t_0)$, hence
$(\bar{q},\hat{Q})(t_0)=(q_1,Q)(t_0)$ and
$(\bar{q},\hat{Q})(t_2)=(q_1,Q)(t_2)$. From the minimality of
$(q_1,Q)$, it follows that
\[\int_{t_0}^{t_2}L_n\left(q_1,Q,\dot{q}_1,\dot{Q}\right)dt\leq\int_{t_0}^{t_2}L_n\left(\bar{q},\hat{Q},\dot{\bar{q}},\dot{\hat{Q}}\right)dt.\]
Since $t_2-t_0\geq 2\tau/\delta'_n$, we let
$\bar{\Omega}=[t_1-\frac{\tau}{\delta'_n}+1,t_1+\frac{\tau}{\delta'_n}-1]$,
then (\ref{separax}) and a direct calculation implies
\begin{align*}
\int_{t_0}^{t_2}(1-\cos(\bar{q}(t)))dt&=\int_{\bar{\Omega}}
1-\cos(\hat{q}_1(t-t_1))dt\\
&\ \ \ +\int_{[t_0,t_2]\backslash{\bar{\Omega}}}(1-\cos(\bar{q}(t)))dt,\\
&\preceq\frac{1}{\sqrt{\delta_n}}\int_0^{2\pi}\sqrt{1-\cos(\hat{q}_1)}d\hat{q}_1,\\
&\preceq\frac{1}{\sqrt{\delta_n}},
\end{align*}
hence,
\begin{equation}\label{fu}
\int_{t_0}^{t_2}(1-\cos(\bar{q}))dt\preceq
\frac{1}{\sqrt{\delta_n}}.
\end{equation}
Moreover, let $\bar{\omega}=(\bar{\omega}_2,\ldots,\bar{\omega}_d)$,
in terms of the definition of $L_n$ as (\ref{L}), we have
\begin{align*}
\int_{t_0}^{t_2}L_n&\left(\bar{q},\hat{Q},\dot{\bar{q}},\dot{\hat{Q}}\right)dt\\
&\leq\int_{t_0}^{t_2}\frac{1}{2}\sum^d_{i=3}\frac{|\bar{\omega}_i|^2}{|l_{ni}|^2}+\frac{|\bar{\omega}_2|^2}{2|k'_n|^2}\\
&\ \ \ +\frac{1}{|k_n|^2}\left(\frac{1}{2}|\dot{\bar{q}}|^2+\delta
_n(1-\cos(\bar{q}))\right)dt
+\frac{1}{|k_n|^2}\mu_n\delta_n\int^{t_2}_{t_0}(1-\cos(\bar{q}))dt,\\
\end{align*}hence, from (\ref{fu}), it follows that
\begin{align*}
\int_{t_0}^{t_2}L_n&\left(\bar{q},\hat{Q},\dot{\bar{q}},\dot{\hat{Q}}\right)dt\\
&\leq\int_{t_0}^{t_2}\frac{1}{2}\sum^d_{i=3}\frac{|\bar{\omega}_i|^2}{|l_{ni}|^2}+\frac{|\bar{\omega}_2|^2}{2|k'_n|^2}\\
&\ \ \ +\frac{1}{|k_n|^2}\left(\frac{1}{2}|\dot{\bar{q}}|^2+\delta
_n(1-\cos(\bar{q}))\right)dt
+\frac{C}{|k_n|^2}\mu_n\sqrt{\delta_n}.\\
\end{align*}Since $\mu_n$ is small enough for $n$ large enough, we
have\begin{align*}
\int_{t_0}^{t_2}L_n&\left(q_1,Q,\dot{q}_1,\dot{Q}\right)dt\\
&\geq\int_{t_0}^{t_2}\frac{1}{2}\sum^d_{i=3}\frac{|\dot{q}_i|^2}{|l_{ni}|^2}+\frac{|\dot{q}_2|^2}{2|k'_n|^2}dt
+\frac{1}{2}\frac{1}{|k_n|^2}\int_{t_0}^{t_2}(|\dot{q}_1|^2+\delta
_n(1-\cos(q_1)))dt.
\end{align*} It is easy to see that
\[\frac{1}{2}|\dot{q}_i|^2=\frac{1}{2}|\bar{\omega}_i|^2+\langle \bar{\omega}_i,
\dot{q}_i-\bar{\omega}_i\rangle+\frac{1}{2}|\dot{q}_i-\bar{\omega}_i|^2,\quad
\text{for}\quad i=2,\ldots,d,\] and the mean value of
$\dot{q}_i-\bar{\omega}_i$ vanishes on $[t_0,t_2]$ for
$i=2,\ldots,d$. Hence,
\begin{align*}
\int_{t_0}^{t_2}L_n&\left(q_1,Q,\dot{q}_1,\dot{Q}\right)dt\\
&\geq
\int_{t_1}^{t_3}\frac{1}{2}\sum^d_{i=3}\frac{|\bar{\omega}_i|^2}{|l_{ni}|^2}+\frac{|\bar{\omega}_2|^2}{2|k'_n|^2}dt
+\frac{1}{2}\frac{1}{|k_n|^2}\int_{t_0}^{t_2}(|\dot{q}_1|^2+\delta
_n(1-\cos(q_1)))dt.
\end{align*}
 From the deduction above, we have
\[\frac{1}{2}\int_{t_0}^{t_2}(|\dot{q}_1|^2+\delta
_n(1-\cos(q_1)))dt\leq
\int_{t_0}^{t_2}\left(\frac{1}{2}|\dot{\bar{q}}|^2+\delta
_n(1-\cos(\bar{q}))\right)dt+C\mu_n\sqrt{\delta_n}.\] From the
definition of $\bar{q}$ and (\ref{separax}), a direct calculation
implies
\begin{align*}
\int_{t_0}^{t_2}&\left(\frac{1}{2}|\dot{\bar{q}}|^2+\delta
_n(1-\cos(\bar{q}))\right)dt\\
&=\int_{\bar{\Omega}}\left(\frac{1}{2}|\dot{\hat{q}}_1(t-t_2)|^2+\delta
_n(1-\cos(\hat{q}_1(t-t_2))\right)dt\\
&\ \ \ +\int_{[t_0,t_2]\backslash\bar{\Omega}}(1-\cos(\bar{q}(t)))dt,\\
&\preceq\sqrt{\delta_n}\int_0^{2\pi}\sqrt{1-\cos(\hat{q}_1)}d\hat{q}_1,\\
&\preceq\sqrt{\delta_n},
\end{align*}
hence,
\begin{equation}\label{0}
\frac{1}{2}\int_{t_0}^{t_2}(|\dot{q_1}|^2+\delta
_n(1-\cos(q_1)))dt\preceq\sqrt{\delta_n}.
\end{equation}
From Euler-Lagrange equation on $(q_1, Q)$, it follows that
\[|\ddot{Q}|\leq\mu_n\delta_n(1-\cos(q_1)).\]
Hence,
\[|\dot{Q}(t)-\dot{Q}(t')|\leq \mu_n\delta_n\int_{t'}^t(1-\cos(q_1))dt,\quad\text{for\ any}\  t,\ t'\in [t_0,t_2].\]
By (\ref{0}), we have
\[|\dot{Q}(t)-\dot{Q}(t')|\preceq\mu_n\sqrt{\delta_n},\quad\text{for\ any}\  t,\ t'\in [t_0,t_2].\]
Moreover,
\[|\dot{Q}(t)-\bar{\omega}|\preceq\mu_n\sqrt{\delta_n},\quad\text{for\ any}\  t\in [t_0,t_2],\]
which yields
\[\sup_\Omega|Q(t)-(\hat{Q}(t)+Q(t_1)-\hat{Q}(t_1))|\preceq\frac{|k_n|^2\mu_n}{\sqrt{\delta_n}},\]
if we integrate $\dot{Q}(t)$ from $t_1$ to $t$. In particular, we
have
\[\sup_\Omega|q_2(t)-(\hat{q}_2(t)+q_2(t_1)-\hat{q}_2(t_1))|\preceq
\frac{|k_n|^2\mu_n}{\sqrt{\delta_n}}.\]Since
\begin{equation}\label{cos} |\cos(\theta_1)-\cos(\theta_2)|\leq
|\theta_1-\theta_2|,
\end{equation} then the proof of (\ref{1}) is completed if we take $n$ large enough.


Second, we prove (\ref{2}).
By (\ref{qqq}), we have
\[\frac{1}{2}\hat{q}_1\left(-\frac{\tau}{\delta'_n}\right)\leq q_1\left(t_1-\frac{\tau}{\delta'_n}\right)\leq \frac{3}{2}\hat{q}_1\left(-\frac{\tau}{\delta'_n}\right),\]
hence, from (\ref{separax}), it follows that
\[q_1\left(t_1-\frac{\tau}{\delta'_n}\right)\succeq\exp\left(-\frac{|k_n|^{\frac{3\epsilon}{4}}\tau}{\sqrt{\delta_n}}\right).\]
Moreover, the above inequality and (\ref{qq1}) imply that for $n$
large enough and $t\in (t_1-\frac{\tau}{\delta'_n},t_2]$,
\begin{equation*}
|q_1(t)-\hat{q}_1(t-t_1)|\preceq\frac{\Gamma}{\delta_n(1-\cos(\hat{q}_1(t-t_1)))^{\frac{3}{2}}}\preceq\epsilon_n.
\end{equation*}

The case in which $t\in [t_1,t_1+\frac{\tau}{\delta'_n}]$ is
similar, more
precisely,\[|q_1(t)-\hat{q}_1(t-t_1)|\preceq\epsilon_n,\quad\text{for}\quad
t\in \left[t_1,t_1+\frac{\tau}{\delta'_n}\right].\] Therefore, from
(\ref{cos}), it follows that
\[|\cos(q_1(t))-\cos(\hat{q}_1(t-t_1))|\preceq\epsilon_n,\quad\text{for}\quad
t\in\Omega,\]hence, we complete the proof of (\ref{2}).

Third, we prove (\ref{3}). It suffices to prove for
$\Omega'=\left[t_{-1}+\frac{\tau}{\delta'_n},t_2-\frac{\tau}{\delta'_n}\right]$
\[\delta_n\int_{\Omega'}|1-\cos(q_1(t))|dt\preceq\epsilon_n,\] where $t_{-1}$ satisfies $q_1(t_{-1})=-\pi$.

Let $\tilde{q}_1$ be the solution of $A_n$ satisfying boundary
conditions
\begin{equation*}
\begin{cases}
\tilde{q}_1\left(t_{-1}+\frac{\tau}{\delta'_n}\right)=q_1\left(t_{-1}+\frac{\tau}{\delta'_n}\right),\\
\\
\tilde{q}_1\left(t_1-\frac{\tau}{\delta'_n}\right)=q_1\left(t_1-\frac{\tau}{\delta'_n}\right).
\end{cases}
\end{equation*}
From the minimality of $(q_1,Q)$, it follows that
\[\int_{\Omega'}L_n\left(q_1,Q,\dot{q}_1,\dot{Q}\right)dt\leq\int_{\Omega'}L_n\left(\tilde{q}_1,Q,\dot{\tilde{q}}_1,\dot{Q}\right)dt.\]
By (\ref{L}), we have
\begin{align*}
\int_{\Omega'}&\frac{1}{2}|\dot{q}_1|^2+\delta_n(1-\cos(q_1))+\mu_n\delta_n(1-\cos(q_1))\cos (q_2) dt\\
&\leq\int_{\Omega'}\frac{1}{2}|\dot{\tilde{q}}_1|^2+\delta_n(1-\cos(\tilde{q}_1))+\mu_n\delta_n(1-\cos(\tilde{q}_1))\cos
(q_2) dt.
\end{align*}
From the construction of $\tilde{q}_1$, let $\tilde{e}$ be the
energy of $\tilde{q}_1$, a similar argument as the one in Lemma
\ref{energ} implies
\[0<\tilde{e}\leq
\delta_n\exp\left(-\frac{C|k_n|^{\epsilon}}{\sqrt{\delta_n}}\right).\]Based
on the change of variable $dq=\dot{q}dt$, a direct calculation gives
\begin{align*}
\int_{\Omega'}&\frac{1}{2}|\dot{\tilde{q}}_1|^2+\delta_n(1-\cos(\tilde{q}_1))+\mu_n\delta_n(1-\cos(\tilde{q}_1))\cos (q_2) dt\\
&\preceq\exp\left(-\frac{C|k_n|^{\epsilon}}{\sqrt{\delta_n}}\right).\end{align*}
For $n$ large enough, we have
\begin{align*}
\frac{1}{2}&\int_{\Omega'}|\dot{q}_1|^2+\delta_n(1-\cos(q_1))dt\\
&\leq\int_{\Omega'}\frac{1}{2}|\dot{q}_1|^2+\delta_n(1-\cos(q_1))+\mu_n\delta_n(1-\cos(q_1))\cos
(q_2) dt,
\end{align*}
which together with $|\dot{q}_1|\geq 0$ implies
\[\delta_n\int_{\Omega'}|1-\cos(q_1(t))|dt\leq\int_{\Omega'}|\dot{q}_1|^2+\delta_n(1-\cos(q_1))dt\preceq\epsilon_n,\]hence, (\ref{3}) is proved.

Finally, the inequality (\ref{4}) can be obtained by a direct
calculation. More precisely,
\begin{align*}\delta_n\int_{\R\backslash
\Omega}|1-\cos(\hat{q}_1(t-t_1))|dt=&\delta_n\int_{-\infty}^{t_1-\frac{\tau}{\delta'_n}}1-\cos(\hat{q}_1(t-t_1))dt\\
&+\delta_n\int^{\infty}_{t_1+\frac{\tau}{\delta'_n}}1-
\cos(\hat{q}_1(t-t_1))dt.\end{align*} By (\ref{separax}),
$\hat{q}_1(t)=4\arctan\left(\exp(\sqrt{\delta_n}(t-t_1))\right)$, we
have
\begin{align*}
\delta_n\int_{-\infty}^{t_1-\frac{\tau}{\delta'_n}}1-&\cos(\hat{q}_1(t-t_1))dt\\
&\preceq\delta_n\int_{-\infty}^{t_1-\frac{\tau}{\delta'_n}}\arctan^2\left(\exp(\sqrt{\delta_n}(t-t_1))\right)dt,\\
&\preceq\delta_n\int_{-\infty}^{t_1-\frac{\tau}{\delta'_n}}\exp\left(2\sqrt{\delta_n}(t-t_1)\right)dt,\\
&\preceq\sqrt{\delta_n}\exp\left(-\frac{2|k_n|^{\frac{3\epsilon}{4}}\tau}{\sqrt{\delta_n}}\right).
\end{align*}
It is similar to get
\[\delta_n\int^{\infty}_{t_1+\frac{\tau}{\delta'_n}}1-
\cos(\hat{q}_1(t-t_1))dt\preceq\sqrt{\delta_n}\exp\left(-\frac{2|k_n|^{\frac{3\epsilon}{4}}\tau}{\sqrt{\delta_n}}\right).\]
Hence, (\ref{4}) is proved for $n$ large enough,.

So far, we complete the proof of (\ref{1})-(\ref{4}) and also the
proof of Lemma \ref{key}.\End
\begin{Remark}\label{9}
In the proof of (\ref{1}), as a bonus we obtain an important
estimate on $\dot{Q}(t)$. Namely, let $(q_1(t),Q(t))$ be the action
minimizing orbit of $L_n$, then for any $t',t''\in \R$ and $t\in
[t',t'']$ we have
\begin{equation*}\left|\dot{Q}(t)-\frac{Q(t'')-Q(t')}{t''-t'}\right|\leq
C\mu_n\sqrt{\delta_n}.\end{equation*}
\end{Remark}
From Remark \ref{9} and the definition of Melnikov function
(\ref{mel}), it follows that for $n$ large enough,
\begin{equation}\label{22}
|M_n(\omega_2, Q_1,t_1)-M_n(\bar{\omega}_2,Q_1,t_1)|\leq
C\sqrt{\delta_n}\exp\left(-\frac{\lambda}{\sqrt{\delta_n}}\right),
\end{equation}
where $\lambda$ is the same as the one in (\ref{mm}).

\paragraph{Proof of Theorem \ref{maint11}} Based on the minimal property of
the orbits on an invariant Lagrangian torus and its graph property,
passing through each $x\in\T^d$, there is a unique minimal curve
$q(t)$ with rotation vector $\omega$ if the Hamiltonian flow
generated by $H_n$ admits a Lagrangian torus with rotation vector
$\omega$. Hence, it is sufficient to prove the existence of some
point in $\T^d$ where no minimal curve passes through.

Indeed, any minimal curve does not pass through the subspace
$(\pi,0)\times\T^{d-2}$. It implies Lemma \ref{key}. Let
 us assume contrary, namely, there exists $\bar{t}_1$ such that
\[q_1(\bar{t}_1)=\pi,\quad q_2(\bar{t}_1)=0,\] where
$q(t)=(q_1,q_2,\ldots,q_d)(t)$ is a minimal orbit in the universal
covering space $\R^d$. Because of $\omega_1\neq 0$, there exist
$t_0$ and $t_1$ such that
\[q_1(t_0)=0,\quad q_1(t_2)=2\pi.\] Obviously, $t_0<\bar{t}_1<t_2$. From (\ref{ap}), we
have
\[t_2-t_0\geq \frac{2\pi}{C_0}|k_n|^{d-1},\] where $C_0$ is the constant in (\ref{ap}). Let $\tilde{t}_1$ be the last time before $\bar{t}_1$ or the first
time after $\bar{t}_1$ such that
\[|q_2(\tilde{t}_1)-q_2(\bar{t}_1)|=\pi.\]
By (\ref{kn1}), we have
\[|\tilde{t}_1-\bar{t}_1|\sim\frac{1}{|\omega_2|}\leq C_1.\] Without loss of generality, one can assume
$\omega_1>0$ and $\omega_2>0$. Consider a solution $\tilde{q}_1$ of
$A_n$ on $(t_0,\tilde{t}_1)$ and on $(\tilde{t}_1,t_2)$ with
boundary conditions respectively
\begin{equation*}
\begin{cases}
\tilde{q}_1(t_0)=q_1(t_0)=0,\\
\tilde{q}_1(\tilde{t}_1)=q_1(\bar{t}_1)=\pi,
\end{cases}
\quad
\begin{cases}
\tilde{q}_1(\tilde{t}_1)=q_1(\bar{t}_1)=\pi,\\
\tilde{q}_1(t_2)=q_1(t_2)=2\pi.
\end{cases}
\end{equation*}

Since $q=(q_1,q_2,\ldots,q_d)$ is assumed to be a minimal curve,
setting $Q=(q_2,\ldots,q_d)$, we have
\begin{equation}\label{ss}\int_{t_0}^{t_2}L_n(\tilde{q}_1,Q,\dot{\tilde{q}}_1,\dot{Q})dt-\int_{t_0}^{t_2}L_n(q_1,Q,\dot{q}_1,\dot{Q})dt\geq
0.\end{equation}

See Fig.3, where $x_1=(q_1(\bar{t}_1),q_2(\bar{t}_1))=(\pi,0)$,
$x_0=(q_1(t_0),q_2(t_0))=(0,q_2(t_0))$,
$x_2=(q_1(t_2),q_2(t_2))=(2\pi,q_2(t_2))$, $\tilde{x}'_1=(\pi,-\pi)$
and $\tilde{x}''_1=(\pi,\pi)$.

\input{figure3.TpX}

$(\tilde{q}_1(t),q_2(t))$ passes through the point $\tilde{x}'_1$ or
$\tilde{x}''_1$.  Thus, by the construction of $L_n$, we obtain from
(\ref{ss}) that
\begin{equation}\label{AandP}
\begin{split}
\int_{t_0}^{t_2}A_n(\tilde{q}_1,\dot{\tilde{q}}_1)dt&-\int_{t_0}^{t_2}A_n(q_1,\dot{q}_1)dt\\
&\geq
\mu_n\left(\int_{t_0}^{t_2}\tilde{P}_n(q_1,Q)dt-\int_{t_0}^{t_2}\tilde{P}_n(\tilde{q}_1,Q)dt\right),
\end{split}
\end{equation}
where
\[\tilde{P}_n(q_1,Q)=\mu_n\delta_n(1-\cos q_1)\cos q_2.\]
By Lemma \ref{key}, we have
\begin{align*}
\int_{t_0}^{t_2}\tilde{P}_n(\tilde{q}_1,Q)dt&-\int_{t_0}^{t_2}\tilde{P}_n(q_1,Q)dt\\
&\leq
M_n(\bar{\omega}_2,Q(\tilde{t}_1),\tilde{t}_1)-M_n(\bar{\omega}_2,Q(\bar{t}_1),\bar{t}_1)+C\sqrt{\delta_n
}\exp\left(-\frac{\lambda}{\sqrt{\delta_n}}\right),
\end{align*} where $\bar{\omega}=(Q(t_2)-Q(t_0))/(t_2-t_0)$. Here, the approximation from $M_n(\bar{\omega}_2,Q(\tilde{t}_1),\tilde{t}_1)$
to $\int_{t_0}^{t_2}\tilde{P}_n(\tilde{q}_1,Q)dt$ can not be
obtained directly by Lemma \ref{key}, since $(\tilde{q}_1, Q(t))$
may be not minimal. But based on the simplicity of $\tilde{q}_1$, a
much simpler calculation than the one in Lemma \ref{key} implies
that the approximation in Lemma \ref{key} is still verified for the
orbit $(\tilde{q}_1, Q(t))$. We omit it.

By Lemma \ref{10} and (\ref{22}), we have
\begin{equation*}
\int_{t_0}^{t_2}\tilde{P}_n(\tilde{q}_1,Q)dt-\int_{t_0}^{t_2}\tilde{P}_n(q_1,Q)dt\preceq
-\exp\left(-\frac{\lambda}{\sqrt{\delta_n}}\right).
\end{equation*}
From (\ref{AandP}), it follows that
\begin{equation}\label{contr1}
\int_{t_0}^{t_2}A_n(\tilde{q}_1,\dot{\tilde{q}}_1)dt-\int_{t_0}^{t_2}A_n(q_1,\dot{q}_1)dt\succeq\mu_n\exp\left(-\frac{\lambda}{\sqrt{\delta_n}}\right).
\end{equation}

On the other hand, consider a solution $\bar{q}_1$ of $A_n$ on
$(t_0,\bar{t}_1)$ and on $(\bar{t}_1,t_2)$ with boundary conditions
respectively
\begin{equation*}
\begin{cases}
\bar{q}_1(t_0)=q_1(t_0)=0,\\
\bar{q}_1(\bar{t}_1)=q_1(\bar{t}_1)=\pi,
\end{cases}
\quad
\begin{cases}
\bar{q}_1(\bar{t}_1)=q_1(\bar{t}_1)=\pi,\\
\bar{q}_1(t_2)=q_1(t_2)=2\pi.
\end{cases}
\end{equation*}
We have
\[\int_{t_0}^{t_2}A_n(q_1,\dot{q}_1)dt\geq\int_{t_0}^{t_2}A_n(\bar{q}_1,\dot{\bar{q}}_1)dt.\]
Since $|\tilde{t}_1-\bar{t}_1|\leq C_1$, An argument similar as the
one to Fig. 2 implies
\[\int_{t_0}^{t_2}A_n(\tilde{q}_1,\dot{\tilde{q}}_1)dt-\int_{t_0}^{t_2}A_n(\bar{q}_1,\dot{\bar{q}}_1)dt\preceq|k_n|^{d-1}
\delta_n\exp\left(-\frac{C|k_n|^{\epsilon}}{\sqrt{\delta_n}}\right).\]
From (\ref{mu}),
\[\mu_n\sim
\exp\left(-|k_n|^{\frac{d}{2}-\frac{1}{2}+\frac{\epsilon}{3}}\right),\]
it follows that
\[\int_{t_0}^{t_2}A_n(\tilde{q}_1,\dot{\tilde{q}}_1)dt-\int_{t_0}^{t_2}A_n(\bar{q}_1,\dot{\bar{q}}_1)dt
\preceq\frac{\mu_n}{|k_n|}\exp\left(-\frac{\lambda}{\sqrt{\delta_n}}\right).\]
Hence, for any $\bar{t}_1\in (t_0,t_2)$, we can find $\tilde{t}_1$
such that
\begin{align*}
\int_{t_0}^{t_2}A_n(\tilde{q}_1,\dot{\tilde{q}}_1)dt-\int_{t_0}^{t_2}A_n(q_1,\dot{q}_1)dt&\leq
\int_{t_0}^{t_2}A_n(\tilde{q}_1,\dot{\tilde{q}}_1)dt-\int_{t_0}^{t_2}A_n(\bar{q}_1,\dot{\bar{q}}_1)dt,\\
&\preceq\frac{\mu_n}{|k_n|}\exp\left(-\frac{\lambda}{\sqrt{\delta_n}}\right),
\end{align*}
which is contradicted by (\ref{contr1}) for large $n$. this
completes the proof of Theorem \ref{maint11}.  \End

\section{\sc Destruction of all Lagrangian tori}
In this section, we are concerned with Lagrangian tori as follow:
\begin{Definition}\label{dd1}
$\mathcal {T}^d$ is called $d$ dimensional Lagrangian torus with the
rotation vector $\omega$ if
\begin{itemize}
\item $\mathcal {T}^d$ is a Lagrangian submanifold;
\item $\mathcal {T}^d$ is invariant for the Hamiltonian
flow $\Phi_H^t$ generated by $H$.
\end{itemize}
\end{Definition}
 For positive definite Hamiltonian
systems, if a Lagrangian torus is invariant under the Hamiltonian
flow, it is then the graph over $\T^d$ (see \cite{BP}).

In the following subsections, we consider the destruction of all
Lagrangian tori of symplectic twist maps. Based on the
correspondence between symplectic twist maps and Hamiltonian
systems, it can be achieved to destruct all Lagrangian tori of
Hamiltonian systems.
\subsection{A toy model} To show the basic ideas, we are beginning
with a toy model whose generating function is as follow:
\begin{equation}\label{hall}
h_n(x,x')=h_0(x,x')-\frac{5}{4n^2}\sin(nx')-\frac{1}{16 n^2}\cos(2
nx'),
\end{equation}
where $h_0(x,x')=\frac{1}{2}(x-x')^2$. Let $f_n(x,y)=(x',y')$ be the
exact area-preserving twist map generated by (\ref{hall}), then
\begin{equation*}
\begin{cases}
y=-\partial_1 h_n(x,x')=x'-x,\\
y'=\partial_2 h_n(x,x')=x'-x-\frac{5}{4 n}\cos( nx')+\frac{1}{8
n}\sin(2 nx').
\end{cases}
\end{equation*}
We set $\phi_n(x)=-\frac{5}{4 n}\cos( nx)+\frac{1}{8 n}\sin(2 nx)$,
then
\begin{equation}
f_n(x,y)=(x+y,y+\phi_n(x+y)).
\end{equation}
In \cite{H1}, Herman found a criterion of total destruction of
invariant circles. By Birkhoff graph theorem (see \cite{H2}), if
$f_n$ admits an invariant circle, then the invariant circle is a
Lipschitz graph. We denote the graph by $\psi_n$, then it follows
from \cite{H3} that
\begin{equation*}
\psi_n\circ g_n-\psi_n=\phi_n\circ g_n,
\end{equation*}
where $g_n=\text{Id}+\psi_n$. This is equivalent to
\begin{equation}\label{gg}
\frac{1}{2}(g_n+g_n^{-1})=\text{Id}+\frac{1}{2}\phi_n.
\end{equation}
Let $\mathfrak{D}_n$ be the set of differentiable points of $g_n$,
then $\mathfrak{D}_n$ has full Lebesgue measure on $\R$ since $g_n$
is a Lipschitz function. For $x\in \mathfrak{D}_n$, we differentiae
(\ref{gg}),
\begin{equation*}
\frac{1}{2}(Dg_n(x)+(Dg_n)^{-1}(g_n^{-1}(x)))=1+\frac{1}{2}D\phi_n(x).
\end{equation*}
Let $G_n=||Dg_n||_{L^\infty}$. It is easy to see that for $\eps>0$,
there exists $\tilde{x}\in \mathfrak{D}_n$ such that
$Dg_n(\tilde{x})\geq G_n-\eps$. Let $M_n=\max D\phi_n$, we have
\begin{equation*}
\frac{1}{2}\left(G_n+\frac{1}{G_n}-\eps\right)\leq 1+\frac{1}{2}M_n.
\end{equation*}
Since $\eps>0$ is arbitrarily small, then
\begin{equation*}
\frac{1}{2}\left(G_n+\frac{1}{G_n}\right)\leq 1+\frac{1}{2}M_n.
\end{equation*}
Hence,
\begin{equation}\label{mmm}
G_n\leq
1+\frac{1}{2}M_n+\left(M_n+\frac{1}{4}{M_n}^2\right)^{\frac{1}{2}}.
\end{equation}
Obviously, for $x\in \mathfrak{D}_n$, we have
\begin{equation*}
\frac{1}{G_n}\leq 1+\frac{1}{2}D\phi_n(x).
\end{equation*}
Let $m_n=\min D\phi_n$, then we have
\begin{equation*}
\frac{1}{G_n}\leq 1+\frac{1}{2}m_n,
\end{equation*}
which together with (\ref{mmm}) implies that
\begin{equation*}
\frac{1}{1+\frac{1}{2}m_n}\leq
1+\frac{1}{2}M_n+\left(M_n+\frac{1}{4}{M_n}^2\right)^{\frac{1}{2}}.
\end{equation*}
Therefore, it is sufficient for total destruction of invariant
circles to construct $\phi_n(x)$ such that
\begin{equation}\label{cri}
\frac{1}{1+\frac{1}{2}\min D\phi_n}> 1+\frac{1}{2}\max
D\phi_n+\left(\max D\phi_n+\frac{1}{4}(\max
D\phi_n)^2\right)^{\frac{1}{2}}.
\end{equation}
In our construction,
\begin{equation}
D\phi_n(x)=\frac{5}{4}\sin(nx)+\frac{1}{4}\cos(2 nx).
\end{equation}
A simple calculation implies
\begin{equation*}
\left\{\begin{array}{ll}
\hspace{-0.4em}\min D\phi_n(x)=-\frac{3}{2},&\text{attained at}\ x=\frac{3}{2n}+\frac{2\pi k}{n},\\
\hspace{-0.4em}\max D\phi_n(x)=1,&\text{attained at}\ x=\frac{\pi}{2n}+\frac{2\pi k}{n},\\
\end{array}\right.
\end{equation*}
where $k\in\Z$. Hence, (\ref{cri}) holds. Moreover, the exact
area-preserving twist map generated by (\ref{hall}) admits no
invariant circles.

By interpolation inequality (\cite{H1}), for a small positive
constant $\delta$, we have
\[||\phi_n||_{C^{1-\delta}}\leq 2||\phi_n||_{C^0}^\delta||D\phi_n||_{C^0}^{1-\delta}.\]
From the construction of $\phi_n$, it follows that
$||\phi_n||_{C^0}\rightarrow 0$, as $n\rightarrow \infty$ and
$||D\phi_n||_{C^0}$ is bounded. Hence,
\begin{equation*}
||\phi_n||_{C^{1-\delta}}\rightarrow 0\quad\text{as}\quad
n\rightarrow\infty,
\end{equation*}
which implies that
\begin{equation*}
||h_n-h_0||_{C^{2-\delta}}\rightarrow 0\quad\text{as}\quad
n\rightarrow\infty.
\end{equation*}
\subsection{$C^\infty$ case}
 In \cite{H3}, Herman extended
the criterion (\ref{cri}) to multi-degrees of freedom. More
precisely, for exact symplectic twist map of $\text{T}^*\T^d$ of the
following form
\begin{equation}
f(x,y)=(x+y,y+d\Psi(x+y)),
\end{equation}
where $\Psi\in C^r(\T^d,\R)$, $r\geq 2$ and
\[d\Psi=\left(\frac{\partial\Psi}{\partial x_1},\cdots,\frac{\partial\Psi}{\partial x_d}\right).\]
Let $E(x)$ be the derivative matrix of $d\Psi$ and
$T(x)=\frac{1}{d}\text{tr}E(x)$, where $\text{tr}$ denotes the trace
of $E(x)$. From a similar argument as the deduction of (\ref{cri}),
it follows that it is sufficient for total destruction of Lagrangian
tori to construct $T(x)$ such that
\begin{equation}\label{cri_multi}
\frac{1}{1+\frac{1}{2}\min T(x)}> 1+\frac{1}{2}\max T(x)+\left(\max
T(x)+\frac{1}{4}(\max T(x))^2\right)^{\frac{1}{2}}.
\end{equation}
Moreover, for $T(x)\rightarrow 0$, (\ref{cri_multi}) implies
\begin{equation}\label{scri}
-\frac{1}{2}\min T(x)>\sqrt{\max T(x)}+O(\max T(x)).
\end{equation}
Herman constructed a sequence $\{\Psi_n\}_{n\in\N}$ that satisfies
(\ref{scri}). It is easy to see $T_n(x)=\frac{1}{d}\Delta\Psi_n$
where $\Delta$ denotes the Laplacian. Since $T_n(x)$ is
$2\pi$-periodic, it is enough to construct it on $[-\pi,\pi]^d$.
More precisely,
\begin{equation*}
T_n(x)=\left\{\begin{array}{ll}
\hspace{-0.4em}T_n^+(x),&x\in [0,\pi]^d,\\
\hspace{-0.4em}T_n^-(x),&x\in [-\pi,0]^d,\\
\hspace{-0.4em}0,&\text{others}.\\
\end{array}\right.
\end{equation*}
where $T_n(x)$ is $C^\infty$ function, $T_n^+(x)$ and $T_n^-(x)$
have the following forms respectively.

$T_n^+(x)$ satisfies:
\begin{equation*}
\left\{\begin{array}{l}
\hspace{-0.4em}\text{supp}\,T_n^+(x)\subset [0,\pi]^d,\\
\hspace{-0.4em}\max T_n^+(x)=\frac{1}{n},\\
\end{array}\right.
\end{equation*}

$T_n^-(x)$ satisfies:
\begin{equation*}
\left\{\begin{array}{l}
\hspace{-0.4em}\text{supp}\,T_n^-(x)=B_{R_n}(x_0),\\
\hspace{-0.4em}\max T_n^-(x)=\frac{1}{\sqrt{n}},\\
\hspace{-0.4em}R_n\sim
\left(\frac{1}{\sqrt{n}}\right)^\frac{1}{d},\\
\hspace{-0.4em}x_0=\left(-\frac{\pi}{2},\cdots,-\frac{\pi}{2}\right),\\
\end{array}\right.
\end{equation*}
where $f\sim g$ means that $\frac{1}{C}g<f<Cg$ for a constant $C>1$.
Hence, we obtain a sequence of $\{T_n(x)\}_{n\in \N}$ with bounded
$C^d$ norms and satisfying
\[\int_{\T^d}T_n(x)dx=0.\]
Form interpolation inequality, it follows that $T_n(x)\rightarrow 0$
as $n\rightarrow \infty$ in the $C^{d-\delta}$ topology for any
$\delta>0$.

Let $\Psi_n$ be the unique function in $C^\infty(\T^d,\R)$ such that
\[\int_{\T^d}\Psi_n(x)dx=0\quad\text{and}\quad\frac{1}{d}\Delta\Psi_n(x)=T_n(x).\]
By Schauder estimates one knows that for any $\delta>0$,
$\Psi_n(x)\rightarrow 0$ as $n\rightarrow\infty$ in the
$C^{d+2-\delta}$ topology. From the construction of $T_n(x)$, it is
easy to see that (\ref{scri}) is verified.

Above all, we have the following theorem
\begin{Theorem}
All Lagrangian tori of an integrable symplectic twist map with
$d\geq 1$ degrees of freedom can be destructed by $C^\infty$
perturbations of the generating function and the perturbations are
arbitrarily small in the $C^{d+2-\delta}$  topology for a small
given constant $\delta>0$.
\end{Theorem}
Based on the correspondence between symplectic twist maps and
Hamiltonian systems, we have the following corollary.
\begin{Corollary}
All Lagrangian tori of an integrable Hamiltonian system with $d\geq
2$ degrees of freedom can be destructed by $C^\infty$ perturbations
which are arbitrarily small in the $C^{d+1-\delta}$ topology  for a
small given constant $\delta>0$.
\end{Corollary}
\subsection{An approximation lemma} In this subsection, we will
prove a lemma on $C^\infty$ functions approximated by trigonometric
polynomials. First of all, we need some notations. Define
\[C^\infty_{2\pi}(\R^d,\R):=\left\{f:\R^d\rightarrow\R|f\in C^\infty(\R^d,\R)\ \text{and}\ 2\pi-\text{periodic in}\ x_1,\ldots,x_d\right\}.\]
 Let $f(x)\in C^\infty_{2\pi}(\R^d,\R)$. The $m$-th
 Fej\'{e}r-polynomial of $f$ with respect to $x_j$ is given by
 \begin{equation}\label{B1}
F_m^{[j]}(f)(x):=\frac{1}{m\pi}\int_{-\pi/2}^{\pi/2}f(x+2e_j)\left(\frac{\sin(mt)}{\sin
t}\right)^2dt,
 \end{equation}
where $x\in\R^d$, $m\in\N$, $j\in\{1,\ldots,d\}$ and $e_j$ is the
$j$-th vector of the canonical basis of $\R^d$. $F_m^{[j]}(f)(x)$ is
a trigonometric polynomial in $x_j$ of degree at most $m-1$. By
\cite{Z},
\[\frac{1}{m\pi}\int_{-\pi/2}^{\pi/2}\left(\frac{\sin(mt)}{\sin
t}\right)^2dt=1,\] hence, from (\ref{B1}), we have
\[||F_m^{[j]}(f)||_{C^0}\leq ||f||_{C^0}.\]
We denote
\[P_m^{[j]}(f):=2F_{2m}^{[j]}(f)-F_m^{[j]}(f).\] It is easy to see
that $P_m^{[j]}(f)$ is a trigonometric polynomial in $x_j$ of degree
at most $2m-1$. Moreover,
 \begin{equation}\label{B2}
||P_m^{[j]}(f)||_{C^0}\leq 3||f||_{C^0},
 \end{equation}
 \begin{equation}\label{B3}
P_m^{[j]}(af+bg)=aP_m^{[j]}(f)+bP_m^{[j]}(g),
 \end{equation}
where $a,b\in\R$ and $f,g\in C^\infty_{2\pi}(\R^d,\R)$. For
$k\in\{1,\ldots,d\}$, $j_1,\ldots,j_k\in\{1,\ldots,d\}$ with
$j_p\neq j_q$ for $p\neq q$. Let $m_1,\ldots,m_k\in\N$ and $f\in
C^\infty_{2\pi}(\R^d,\R)$, we define
 \begin{equation}\label{B4}
P_{m_1,\ldots,m_k}^{[j_1,\ldots,j_k]}(f):=P_{m_1}^{[j_1]}\left(P_{m_2}^{[j_2]}\left(\cdots\left(P_{m_k}^{[j_k]}(f)\right)\cdots\right)\right).
 \end{equation}
It is easy to see that for all $l\in \{1,\ldots,k\}$,
$P_{m_1,\ldots,m_k}^{[j_1,\ldots,j_k]}(f)$ are trigonometric
polynomials in $x_l$ of degree at most $2m_l-1$. We have the
following lemma.
\begin{Lemma}\label{appp}
Let $f\in C^\infty_{2\pi}(\R^d,\R)$, $r_1,\ldots,r_d\in\N$,
$m_1,\ldots,m_d\in\N$, then we have
\begin{equation}\label{B8}
\|f-P_{m_1,\ldots,m_d}^{[1,\ldots,d]}(f)\|_{C^0}\leq
C_d\sum_{j=1}^{d}\frac{1}{{m_j}^{r_j}}\left\|\frac{\partial
^{r_j}f}{\partial {x_j}^{r_j}}\right\|_{C^0},
\end{equation}
where $C_d$ is a constant only depending on $d$.
\end{Lemma}
\Proof From \cite{A}, we will prove Lemma \ref{appp} by induction.
The case $d=1$ is covered by Jackson's approximation theorem. More
precisely, for $f\in C^\infty_{2\pi}(\R,\R)$, $m,r\in\N$, we have
\begin{equation}\label{B5}
\|f-P_m^{[1]}(f)\|_{C^0}\leq
C_1\frac{1}{{m}^{r}}\left\|\frac{\partial ^{r}f}{\partial
{x}^{r}}\right\|_{C^0}.
\end{equation}
Let the assertion be true for some $d\in\N$. We verify it for $d+1$.
Consider the functions $f(x_1,\cdot)$ with $x_1$ as a real
parameter. Then by the assertion for $d$, we have
\[\|f(x_1,\cdot)-P_{m_2,\ldots,m_{d+1}}^{[2,\ldots,d+1]}(f)(x_1,\cdot)\|_{C^0}\leq
C_d\sum_{j=2}^{d+1}\frac{1}{{m_j}^{r_j}}\left\|\frac{\partial
^{r_j}f}{\partial {x_j}^{r_j}}\right\|_{C^0},\] hence,
\begin{equation}\label{B6}
\|f-P_{m_2,\ldots,m_{d+1}}^{[2,\ldots,d+1]}(f)\|_{C^0}\leq
C_d\sum_{j=2}^{d+1}\frac{1}{{m_j}^{r_j}}\left\|\frac{\partial
^{r_j}f}{\partial {x_j}^{r_j}}\right\|_{C^0}.
\end{equation}
Let $\hat{x}_j\in\R^d$ denote the vector $x\in\R^{d+1}$ without its
$j$-th entry. For the functions $f(\cdot,\hat{x}_1)$, from
(\ref{B5}), it follows that
\[\|f(\cdot,\hat{x}_1)-P_{m_1}^{[1]}(f)(\cdot,\hat{x}_1)\|_{C^0}\leq C_1\frac{1}{{m_1}^{r_1}}\left\|\frac{\partial ^{r_1}f}{\partial
{x_1}^{r_1}}\right\|_{C^0},\]hence,
\begin{equation}\label{B7}
\|f-P_{m_1}^{[1]}(f)\|_{C^0}\leq
C_1\frac{1}{{m_1}^{r_1}}\left\|\frac{\partial ^{r_1}f}{\partial
{x_1}^{r_1}}\right\|_{C^0}.
\end{equation}
By (\ref{B2}), (\ref{B3})),(\ref{B4}) and (\ref{B6}), we have
\begin{align*}
\left\|P_{m_1}^{[1]}(f)-P_{m_1,\ldots,m_{d+1}}^{[1,\ldots,d+1]}(f)\right\|_{C^0}&=\left\|P_{m_1}^{[1]}(f)-P_{m_1}^{[1]}\left(P_{m_2,\ldots,m_{d+1}}^{
[2,\ldots,d+1]}(f)\right)\right\|_{C^0},\\
&=\left\|P_{m_1}^{[1]}\left(f-P_{m_1}^{[1]}P_{m_2,\ldots,m_{d+1}}^{
[2,\ldots,d+1]}(f)\right)\right\|_{C^0},\\
&\leq 3C_d\sum_{j=2}^{d+1}\frac{1}{{m_j}^{r_j}}\left\|\frac{\partial
^{r_j}f}{\partial {x_j}^{r_j}}\right\|_{C^0},
\end{align*}
which together with (\ref{B7}) implies that
\begin{align*}
\left\|f-P_{m_1,\ldots,m_{d+1}}^{[1,\ldots,d+1]}(f)\right\|_{C^0}&\leq\left\|f-P_{m_1}^{[1]}(f)\right\|_{C^0}+\left\|P_{m_1}^{[1]}(f)-P_{m_1}^{[1]}\left(P_{m_2,\ldots,m_{d+1}}^{
[2,\ldots,d+1]}(f)\right)\right\|_{C^0},\\
&\leq C_1\frac{1}{{m_1}^{r_1}}\left\|\frac{\partial
^{r_1}f}{\partial
{x_1}^{r_1}}\right\|_{C^0}+3C_d\sum_{j=2}^{d+1}\frac{1}{{m_j}^{r_j}}\left\|\frac{\partial
^{r_j}f}{\partial {x_j}^{r_j}}\right\|_{C^0},\\
&\leq
C_{d+1}\sum_{j=1}^{d+1}\frac{1}{{m_j}^{r_j}}\left\|\frac{\partial
^{r_j}f}{\partial {x_j}^{r_j}}\right\|_{C^0}.
\end{align*}
This finishes the proof. \End

Obviously, there exist $m_{\bar{j}}, r_{\bar{j}}$ such that
\[\frac{1}{{m_{\bar{j}}}^{r_{\bar{j}}}}\left\|\frac{\partial
^{r_{\bar{j}}}f}{\partial
{x_{\bar{j}}}^{r_{\bar{j}}}}\right\|_{C^0}=\max_{1\leq j \leq
d}\left\{\frac{1}{{m_j}^{r_j}}\left\|\frac{\partial
^{r_j}f}{\partial {x_j}^{r_j}}\right\|_{C^0}\right\}.\] Hence, we
have
\begin{align*}
\|f-P_{m_1,\ldots,m_d}^{[1,\ldots,d]}(f)\|_{C^0}&\leq
dC_d\frac{1}{{m_{\bar{j}}}^{r_{\bar{j}}}}\left\|\frac{\partial
^{r_{\bar{j}}}f}{\partial
{x_{\bar{j}}}^{r_{\bar{j}}}}\right\|_{C^0},\\
&\leq
C'_d\frac{1}{{m_{\bar{j}}}^{r_{\bar{j}}}}\|f\|_{C^{r_{\bar{j}}}}.
\end{align*}
For the simplicity of notations, we denote
\[p_N(x)=P_{m_1,\ldots,m_d}^{[1,\ldots,d]}(f)(x),\]
where $x=(x_1,\ldots,x_d)$ and $N=2m_{\bar{j}}-1$. Moreover, we
denote $k:=r_{\bar{j}}$, then
\begin{equation}\label{B9}
\|f(x)-p_N(x)\|_{C^0}\leq A_{dk}N^{-k}\|f(x)\|_{C^k},
\end{equation}
where $A_{dk}$ is a constant depending on $d$ and $k$.
\subsection{$C^\omega$ case} Similar to Herman's construction,  we
consider $C^\infty$ function $\tilde{T}_n(x)$ as follow:
\begin{equation*}
\tilde{T}_n(x)=\left\{\begin{array}{ll}
\hspace{-0.4em}\tilde{T}_n^+(x),&x\in [0,\pi]^d,\\
\hspace{-0.4em}\tilde{T}_n^-(x),&x\in [-\pi,0]^d,\\
\hspace{-0.4em}0,&\text{others}.\\
\end{array}\right.
\end{equation*}

$\tilde{T}_n^+(x)$ satisfies:
\begin{equation*}
\left\{\begin{array}{l}
\hspace{-0.4em}\text{supp}\,\tilde{T}_n^+(x)\subset [0,\pi]^d,\\
\hspace{-0.4em}\max \tilde{T}_n^+(x)\sim 1,\\
\end{array}\right.
\end{equation*}

$\tilde{T}_n^-(x)$ satisfies:
\begin{equation*}
\left\{\begin{array}{l}
\hspace{-0.4em}\text{supp}\,\tilde{T}_n^-(x)=B_{R_n}(x_0),\\
\hspace{-0.4em}\max \tilde{T}_n^-(x)=n,\\
\hspace{-0.4em}R_n\sim
\left(\frac{1}{n}\right)^\frac{1}{d},\\
\hspace{-0.4em}x_0=\left(-\frac{\pi}{2},\cdots,-\frac{\pi}{2}\right).\\
\end{array}\right.
\end{equation*}

By Lemma \ref{appp}, there exists a trigonometric polynomial
$p_N(x_1,\cdots,x_d)$ in $x_l$ $(1\leq l \leq d)$ of degree at most
$N$ such that
\[\left\|\tilde{T}_n(x_1,\cdots,x_d)-p_N(x_1,\cdots,x_d)\right\|_{C^0}\leq A_{dk}N^{-k}\left\|\tilde{T}_n(x_1,\cdots,x_d)\right\|_{C^k}.\]
By the construction of $\tilde{T}_n$, we have
\begin{equation}\label{Tck}
||\tilde{T}_n(x_1,\cdots,x_d)||_{C^k}\sim n^{\frac{k}{d}+1}.
\end{equation}
Then, choosing $N$ large enough such that
\begin{equation}\label{N}
A_{dk}N^{-k}||\tilde{T}_n(x_1,\cdots,x_d)||_{C^k}<\sigma\ll 1,
\end{equation}
where $\sigma$ is a small enough positive constant. Hence, we have
\begin{equation}
\left\{\begin{array}{ll}\hspace{-0.4em}\max p_N(x)\sim n,&\text{attained on}\ B_{R_n}(x_0),\\
\hspace{-0.4em}\max p_N(x)\sim 1,&\text{on}\ [-\pi,\pi]^d\backslash\ B_{R_n}(x_0).\\
\end{array}\right.
\end{equation}
By (\ref{Tck}) and (\ref{N}), we have
\begin{equation}\label{nn}
N>\left(\frac{A_{dk}}{\sigma}\right)^{\frac{1}{k}}n^{\frac{1}{d}+\frac{1}{k}}\geq
Cn^{\frac{1}{d}+\frac{1}{k}},
\end{equation}
where $C$ is a constant independent of $n$. We consider the
normalized trigonometric polynomial
\begin{equation}
\tilde{p}_N(x)=\frac{1}{n^{1-\eps}}\frac{p_N(x)}{\max |p_N(x)|},
\end{equation}where $x=(x_1,\ldots,x_d)$. It is easy to see that
\begin{equation*}
\left\{\begin{array}{l}
\hspace{-0.4em}\max \tilde{p}_N(x)\sim \frac{1}{n^{2-\eps}},\\
\hspace{-0.4em}\min \tilde{p}_N(x)\sim -\frac{1}{n^{1-\eps}}.\\
\end{array}\right.
\end{equation*}
Hence, (\ref{scri}) is verified.

Next, we estimate $||\tilde{p}_N(x)||_{C^r}$. By a simple
calculation, we have
\begin{equation}
||p_N(x)||_{C^r}\leq CnN^{r+1}.
\end{equation}
Then,
\begin{align*}
||\tilde{p}_N(x)||_{C^r}&\leq\frac{1}{n^{1-\eps}}\frac{1}{\max
|p_N(x)|} ||p_N(x)||_{C^r},\\
&\leq \frac{1}{n^{2-\eps}}\cdot CnN^{r+1},\\
&=C\frac{1}{n^{1-\eps}}N^{r+1}.
\end{align*}
To achieve $n^{-(1-\eps)}N^{r+1}\rightarrow 0$ as
$n\rightarrow\infty$, it suffices to make
\[\frac{1}{n^{1-\eps}}N^{r+1}\leq \frac{1}{n^\eps}.\] Hence, we have
\[r\leq \log_{N}n^{1-2\eps}-1.\]
From (\ref{nn}), we take
\[N\sim n^{\frac{1}{d}+\frac{1}{k}},\]then, we have
\begin{equation}
\max_N r\leq \left(\frac{1}{d}+\frac{1}{k}\right)^{-1}(1-2\eps)-1.
\end{equation}
Since $k$ can be made large enough, we take $k=\frac{d}{2\eps}$. Let
\[\delta=\frac{4\eps}{1+2\eps}d,\]then, we have
\begin{equation}
\max_N r\leq d-1-\delta.
\end{equation}
It follows that $\tilde{p}_N(x)\rightarrow 0$ as
$n\rightarrow\infty$ in the $C^{d-1-\delta}$ for any $\delta>0$.

Let $\tilde{\Psi}_n$ be the unique function in $C^\omega(\T^d,\R)$
such that
\[\int_{\T^d}\tilde{\Psi}_n(x)dx=0\quad\text{and}\quad\frac{1}{d}\Delta\tilde{\Psi}_n(x)=\tilde{p}_N(x).\]
By Schauder estimates one knows that for any $\delta>0$,
$\tilde{\Psi}_n(x)\rightarrow 0$ as $n\rightarrow\infty$ in the
$C^{d+1-\delta}$ topology. From the construction of
$\tilde{p}_N(x)$, it is easy to see that (\ref{scri}) is verified.
Hence, we have that the symplectic twist maps generated by the
generating function
$\tilde{h}_n(x,x')=\frac{1}{2}(x-x')^2+\tilde{\Psi}_n(x')$ do not
admit any Lagrangian tori for $n$ large enough.

So far, we prove the following theorem
\begin{Theorem}
All Lagrangian tori of an integrable symplectic twist map with
$d\geq 2$ degrees of freedom can be destructed by $C^\omega$
perturbations of the generating function and the perturbations are
arbitrarily small in the $C^{d+1-\delta}$ topology  for a small
given constant $\delta>0$.
\end{Theorem}
Based on the correspondence between symplectic twist maps and
Hamiltonian systems, together with the toy model corresponding to
the case with $d=1$, we have the following corollary.
\begin{Corollary}
All Lagrangian tori of an integrable Hamiltonian system with $d\geq
2$ degrees of freedom can be destructed by $C^\omega$ perturbations
which are arbitrarily small in the $C^{d-\delta}$ topology for a
small given constant $\delta>0$.
\end{Corollary}

 \vspace{2ex}
\noindent\textbf{Acknowledgement} We would like to thank Prof. C.-Q.
Cheng for many helpful discussions. This work is under the support
of Research and Innovation Project for College Graduates of Jiangsu
Province (CXZZ12$\_$0030), the NNSF of China (Grant 11171146) and
PAPD of Jiangsu Province of China.

{\sc Department of Mathematics, Nanjing University, Nanjing 210093,
China.}

 {\it E-mail address:} \texttt{linwang.math@gmail.com}


\addcontentsline{toc}{section}{\sc References}

\begin{thebibliography}{DGOW}
\renewcommand{\baselinestretch}{1.0}
\setlength\itemsep{-1pt}
\small

\bibitem[A]{A} J.Albrecht. {\it On the existence of invariant tori in nearly-integrable Hamiltonian systems with finitely
differentiable perturbations}. Regular and Chaotic Dynamics.
\textbf{12} (2007), 281-320.

\bibitem[Ba]{B} V.Bangert. {\it Mather sets for twist maps and geodesics on
tori}. Dynamics Reported \textbf{1} (1988), 1-45.

\bibitem[Be]{B2} U.Bessi. {\it An analytic counterexample to KAM theorem}.
Ergod. Th. \& Dynam. Sys. \textbf{20} (2000), 317-333.


\bibitem[BP]{BP} M.Bialy and L.Polterovich. {\it Hamiltonian systems,
Lagrangian tori and Birkhoff's theorem}. Math. Ann. \textbf{292}
(1992), 619-627.





\bibitem[Ch]{C} C.-Q.Cheng. {\it Non-existence of KAM torus}. Acta Mathmatica
Sinica. \textbf{27} (2011), 397-404.

\bibitem[Ci]{Chi} B.V.Chirikov. {\it A universal instability of many-dimensional oscillator systems}. Phys. Reports.
\textbf{52} (1979), 264-379.

\bibitem[CW]{CW}
C.-Q.Cheng and L.Wang. {\it Destruction of Lagrangian torus in
positive definite Hamiltonian systems}. Preprint (2012).


\bibitem[F]{Fo} G.Forni. {\it
Analytic destruction of invariant circles}. Ergod.
 Th. \& Dynam. Sys. \textbf{14} (1994), 267-298.

\bibitem[Go]{Go} C.Gol\'{e}. {\it
Optical Hamiltonians and symplectic twist maps}. Physica D:
Nonlinear Phenomena. \textbf{71} (1994), 185-195.

\bibitem[Gr]{Gr} J.M.Greene {\it
A method for determining a stochastic transition}. J.Math. Phys.
\textbf{20} (1979), 1183-1201.

 \bibitem[H1]{H11} M.R.Herman. {\it Sur la conjugation diff$\acute{e}$rentiable des diff$\acute{e}$omorphismes du cercle $\grave{a}$ des
 rotations}. Publ. Math. IHES \textbf{49} (1979), 5-233.

\bibitem[H2]{H1} M.R.Herman. {\it Sur les courbes invariantes par les diff\'{e}omorphismes de
 l'anneau}. Ast\'{e}risque \textbf{103-104} (1983), 1-221.

\bibitem[H3]{H33} M.R.Herman. {\it Sur les courbes invariantes par les
diff$\acute{e}$omorphismes de
 l'anneau}. Ast$\acute{\text{e}}$risque \textbf{144} (1986), 1-243.

\bibitem[H4]{H2} M.R.Herman. {\it In\'{e}galit\'{e}s ``a priori" pour des
tores lagrangiens invariants par des diff\'{e}omrphismes
symplectiques}. Inst. Hautes \'{E}tudes Sci. Publ. Math. \textbf{70}
(1990), 47-101.

\bibitem[H5]{H3} M.R.Herman. {\it Non existence of Lagrangian graphs}.
unpublished preprint (1990).

 \bibitem[Ma1]{M1} J.N.Mather. {\it Existence of quasi periodic orbits for twist homeomorphisms of the
 annulus}. Topology \textbf{21} (1982), 457-467.

 \bibitem[Ma2]{Mm2} J.N.Mather. {\it Non-existence of invariant circles}. Ergod.
 Th. \& Dynam. Sys. \textbf{4} (1984), 301-309.


 \bibitem[Ma3]{M2} J.N.Mather. {\it A criterion for the non-existence of
 invariant circle}. Publ. Math. IHES \textbf{63} (1986), 301-309.

 \bibitem[Ma4]{M3} J.N.Mather. {\it Modulus of continuity for Peierls's
 barrier}. Periodic Solutions of Hamiltonian Systems and Related
 Topics. ed. P.H.Rabinowitz {\it et al}. NATO ASI Series C \textbf{209}.
 Reidel: Dordrecht, (1987), 177-202.

 \bibitem[Ma5]{M4} J.N.Mather. {\it Destruction of invariant circles}. Ergod.
 Th. \& Dynam. Sys. \textbf{8} (1988), 199-214.



\bibitem[Ma6]{M5} J.N.Mather. {\it Action minimizing invariant measures for postive definite Lagrangian
 systems}. Math. Z. \textbf{207} (1991), 169-207.

\bibitem[Mo1]{Mo1} J.Moser. {\it On invariant curves of area-preserving mappings of an annulus}. Nachr. Akad. Wiss. G\"{o}ttingen,
II. Math-Phys. KL. \textbf{1} (1962), 1-20.


\bibitem[Mo2]{Mo2} J.Moser. {\it A rapidly convergent iteration method and non-linear partail differential equations I, II}. Annali della Scuola Normale Superiore di
Pisa. \textbf{20} (1966), 265-315.



\bibitem[Mo3]{Mo3} J.Moser. {\it On the construction of almost periodic solutions for ordinary differential equations}. Proceedings of the International
Conference on Functional Analysis and Related Topics, Tokyo. (1969),
60-67.

 \bibitem[Mo4]{Mo4} J.Moser. {\it Monotone twist mappings and the calculus of variations}. Ergod.
 Th. \& Dynam. Sys. \textbf{6} (1986), 401-413.

\bibitem[MP]{MP} R.S.MacKay and I.C.Percival {\it Converse KAM: theory and practice}. Commun. Math. Phys. \textbf{98} (1985), 469-512.


\bibitem[P]{P} J.P\"{o}schel.  {\it Integrability of Hamiltonian systems on Cantor sets}. Comm. Pure Appl. Math. \textbf{35} (1982),  653-696.

\bibitem[R1]{R1} H.R\"{u}ssman.  {\it On optimal estimates for the solutions of linear partial differential equations of first order with constant
coefficients on the torus}. Lecture Notes in Physics. \textbf{38}
(1975), 598-624.

\bibitem[R2]{R2} H.R\"{u}ssman.  {\it On the existence of invariant curves of twist mappings of an annulus}.
Lecture Notes in Mathematics.  \textbf{1007} (1983),  677-718.


\bibitem[W1]{W1} L.Wang. {\it Variational destruction of invariant
circles. }Discrete and Continuous Dynamical Systems-A. \textbf{32}
(2012), 4429-4443.

\bibitem[W2]{W2} L.Wang. {\it Analytic destruction of invariant circles for exact area-preserving twist
maps}. Preprint (2012).

\bibitem[W3]{W3} L.Wang. {\it Analytic destruction of Lagrangian torus}. Preprint (2012).

\bibitem[W4]{W4} L.Wang. {\it Total destruction of Lagrangian tori}. Preprint (2012).

\bibitem[Z]{Z} A.Zygmund. {\it Trigonometric Series}. Third Edition Volumes I \& II combined, with a foreword by Robert Fefferman. Cambridge
University Press, Cambridge, 2002.
\end{thebibliography}
\end{document}